\input amstex
\input amsppt.sty
\magnification=\magstep1
\hsize=33truecc
\vsize=22.2truecm
\baselineskip=16truept
\NoBlackBoxes
\nologo
\pageno=1
\topmatter
\TagsOnRight

\def\N{\Bbb N}
\def\Z{\Bbb Z}
\def\Q{\Bbb Q}

\def\l{\left}
\def\r{\right}
\def\b{\bigg}
\def\bg{\bigg}
\def\({\b(}
\def\[{\b[}
\def\){\b)}
\def\]{\b]}

\def\t{\text}
\def\f{\frac}
\def\mo{\roman{mod}}

\def\em{\emptyset}
\def\se {\subseteq}

\def\sm{\setminus}

\def\bi{\binom}
\def\eq{\equiv}

\def\ls{\leqslant}
\def\gs{\geqslant}

\def\ve{\varepsilon}
\def\da{\delta}

\def\Proof{\noindent{\it Proof}}
\def\Remark{\noindent{\it Remark}}

\def\Ack{\noindent {\bf Acknowledgments}}
\hbox{Sci. China Math. 58(2015), no.\,7, 1367--1396.}
\medskip
\title On universal sums of polygonal numbers\endtitle
\author Zhi-Wei SUN \endauthor
\affil Department of Mathematics, Nanjing University
     \\Nanjing 210093, People's Republic of China
    \\  zwsun\@nju.edu.cn
    \\ {\tt http://math.nju.edu.cn/$\sim$zwsun}
 \endaffil
\abstract For $m=3,4,\ldots$, the polygonal numbers of order $m$
are given by $p_m(n)=(m-2)\bi n2+n\ (n=0,1,2,\ldots)$. For positive integers $a,b,c$ and $i,j,k\gs3$
with $\max\{i,j,k\}\gs5$, we call the triple $(ap_i,bp_j,cp_k)$ universal if for any $n=0,1,2,\ldots$ there are
nonnegative integers $x,y,z$ such that $n=ap_i(x)+bp_j(y)+cp_k(z)$. We show that there are only 95 candidates
for universal triples (two of which are $(p_4,p_5,p_6)$ and $(p_3,p_4,p_{27})$),
and conjecture that they are indeed universal triples.
For many triples $(ap_i,bp_j,cp_k)$
(including $(p_3,4p_4,p_5),(p_4,p_5,p_6)$ and $(p_4,p_4,p_5)$), we prove that
any nonnegative integer can be
written in the form $ap_i(x)+bp_j(y)+cp_k(z)$ with $x,y,z\in\Z$. We also show some related new results
on ternary quadratic forms, one of which states that any nonnegative integer $n\eq 1\ (\mo\ 6)$ can be written
in the form $x^2+3y^2+24z^2$ with $x,y,z\in\Z$.
In addition, we pose several related conjectures one of which states that for any $m=3,4,\ldots$ each natural number
can be expressed as $p_{m+1}(x_1)+p_{m+2}(x_2)+p_{m+3}(x_3)+r$ with $x_1,x_2,x_3\in\{0,1,2,\ldots\}$
and $r\in\{0,\ldots,m-3\}$.
\endabstract
\keywords Polygonal numbers, ternary quadratic forms, representations of integers, primes\endkeywords
\thanks 2010 {\it Mathematics Subject Classification}.
Primary 11E25; Secondary 11A41, 11B75, 11D85, 11E20, 11P32.
\endthanks
\endtopmatter
\document

\heading{1. Introduction}\endheading

Polygonal numbers are nonnegative integers constructed geometrically from the regular polygons.
For $m=3,4,\ldots$ those {\it $m$-gonal numbers} (or {\it polygonal numbers of order $m$}) are given
by
$$p_m(n):=(m-2)\bi n2+n=\f{(m-2)n^2-(m-4)n}2\
\ (n=0,1,2,\ldots).$$
Clearly,
$$p_m(0)=0,\ p_m(1)=1,\ p_m(2)=m,\ p_m(3)=3m-3,\ p_m(4)=6m-8.$$
Note also that
$$p_3(n)=\f{n(n+1)}2,\ p_4(n)=n^2,\ p_5(n)=\f{3n^2-n}2,\ p_6(n)=2n^2-n.$$
Lagrange's theorem asserts that every $n\in\N=\{0,1,2,\ldots\}$ is the sum of four squares, and Gauss
proved in 1796 a conjecture of Fermat which states that any $n\in\N$ is the sum of three triangular numbers
(this follows from the Gauss-Legendre theorem (see, e.g., [G, pp.\,38-49] or [N96, pp.\,17-23])
which asserts that any positive integer not of the form
$4^k(8l+7)$ with $k,l\in\N$ is the sum of three squares).
Fermat's claim that each $n\in\N$  can be written as the sum of $m$ polygonal numbers of order $m$
was completely proved by Cauchy in 1813 (see [MW, pp.\,54-57] for a proof). Legendre showed that every sufficiently large integer
is the sum of five polygonal numbers of order $m$.
The reader is referred to [N87] and Chapter 1 of [N96, pp.\,3-34] for details.

For $m=3,4$ clearly $p_m(\Z)=\{p_m(x):\ x\in\Z\}$ coincides with $p_m(\N)=\{p_m(n):\ n\in\N\}$.
However, for $m=5,6,\ldots$ we have $p_m(-1)=m-3\in p_m(\Z)\sm p_m(\N)$.
Those $p_m(x)$ with $x\in\Z$ are called {\it generalized $m$-gonal numbers}.
The generalized pentagonal numbers play an important role in combinatorics due to the following celebrated result of
Euler (see, e.g., Berndt [B, p.\,12]):
$$\f1{\sum_{n=0}^{\infty}p(n)q^n}=\prod_{n=1}^\infty(1-q^n)=\sum_{k=-\infty}^\infty (-1)^kq^{p_5(k)}\ \ (|q|<1),$$
where $p(n)$ is the well-known partition function.

As usual, for $A_1,\ldots,A_n\se\Z$  we adopt the notation
$$A_1+\cdots+A_n:=\{a_1+\cdots+a_n:\, a_1\in A_1,\ldots,a_n\in A_n\}.$$
For $a,b,c\in\Z^+=\{1,2,3,\ldots\}$ and $i,j,k\in\{3,4,\ldots\}$, we define
$$N(ap_i,bp_j,cp_k):=ap_i(\N)+bp_j(\N)+cp_k(\N)\tag1.1$$
and
$$Z(ap_i,bp_j,cp_k):=ap_i(\Z)+bp_j(\Z)+cp_k(\Z).\tag1.2$$
If $N(ap_i,bp_j,cp_k)=\N$, then we call the triple $(ap_i,bp_j,cp_k)$ (or the sum $ap_i+bp_j+cp_k$) {\it universal} (over $\N$).
When $Z(ap_i,bp_j,cp_k)=\N$, we say that the sum $ap_i+bp_j+cp_k$ is {\it universal over $\Z$}.

In 1862 Liouville (cf. [D99b, p.\,23]) determined all those universal $(ap_3,bp_3,cp_3)$ (with $1\ls a\ls b\ls c$):
$$\align&(p_3,p_3,p_3),\ (p_3,p_3,2p_3),\ (p_3,p_3,4p_3),\ (p_3,p_3,5p_3),
\\& (p_3,2p_3,2p_3),\ (p_3,2p_3,3p_3),\ (p_3,2p_3,4p_3).
\endalign$$

In 2007 Z. W. Sun [S07] suggested the determination of those universal $(ap_3,bp_3,cp_4)$ and $(ap_3,bp_4,cp_4)$,
and this was  completed via  a series of papers by Sun and his coauthors (cf. [S07], [GPS] and [OS]).
Here is the list of universal triples $(ap_i,bp_j,cp_k)$ with $\{i,j,k\}=\{3,4\}$:
 $$\align&(p_3,p_3,p_4),\ (p_3,p_3,2p_4),\ (p_3,p_3,4p_4),\ (p_3,2p_3,p_4),\ (p_3,2p_3,2p_4),
 \\&(p_3,2p_3,3p_4),\ (p_3,2p_3,4p_4),\ (2p_3,2p_3,p_4),\ (2p_3,4p_3,p_4),\ (2p_3,5p_3,p_4),
 \\&(p_3,3p_3,p_4),\ (p_3,4p_3,p_4),\ (p_3,4p_3,2p_4), \ (p_3,6p_3,p_4),\ (p_3,8p_3,p_4),
 \\&(p_3,p_4,p_4),\ (p_3,p_4,2p_4),\ (p_3,p_4,3p_4),\ (p_3,p_4,4p_4),\ (p_3,p_4,8p_4),
 \\&(p_3,2p_4,2p_4),\ (p_3,2p_4,4p_4),\ (2p_3,p_4,p_4),\ (2p_3,p_4,2p_4),\ (4p_3,p_4,2p_4).
 \endalign$$
 For {\it almost universal} triples $(ap_i,bp_j,cp_k)$ with $\{i,j,k\}=\{3,4\}$, the reader may consult
 [KS] by Kane and Sun, and the papers [CH] and [CO] where the authors proved [KS, Conjecture 1.19] concerning the remaining cases.

In this paper we study universal sums $ap_i+bp_j+cp_k$ with $\max\{i,j,k\}\gs5$. Note that
polygonal numbers of order $m\gs5$ are more sparse than squares and triangular numbers.

Observe that
$$p_6(\N)\se p_6(\Z)=p_3(\Z)=p_3(\N)\tag1.3$$
since
$$p_6(x)=x(2x-1)=p_3(2x-1)\ \t{and}\ p_3(x)=p_6\l(-\f x2\r)=p_6\l(\f{x+1}2\r).$$

It is interesting to
determine all those universal sums $ap_k+bp_k+cp_k$ over $\Z$ with
$a,b,c\in\Z^+=\{1,2,3,\ldots\}$. Note that the case $k=6$ is equivalent to the case $k=3$ which was handled by
Liouville (cf. [D99b, p.\,23]).  That $p_5+p_5+p_5$ is universal over $\Z$ (equivalently,
for any $n\in\N$ we can write $24n+3=x^2+y^2+z^2$ with $x,y,z$ all relatively prime to $3$),
was first realized by R. K. Guy [Gu]. However, as H. Pan told the author, one needs the following identity
$$(1^2+1^2+1^2)^2(x^2+y^2+z^2)=(x-2y-2z)^2+(y-2x-2z)^2+(z-2x-2y)^2\tag1.4$$
(which is a special case of R\'ealis' identity [D99b, p.\,266])
to show that if $x,y,z$ are not all divisible by 3 then
$9(x^2+y^2+z^2)=u^2+v^2+w^2$ for some $u,v,w\in\Z$ with $u,v,w$ not all divisible by 3.

Throughout this paper, for $a,b\in\Z$ we set $[a,b]=\{x\in\Z:\ a\ls x\ls b\}$.

Now we state our first theorem.

\proclaim{Theorem 1.1} {\rm (i)} Suppose that $ap_k+bp_k+cp_k$ is universal over $\Z$,
where $k\in\{4,5,7,8,9,\ldots\}$, $a,b,c\in\Z^+$ and $a\ls b\ls c$.
Then $k$ is equal to $5$ and $(a,b,c)$ is among the following $20$ triples:
$$\align&(1,1,c)\ (c\in[1,10]\sm\{7\}),\ (1,2,c)\ (c\in\{2,3,4,6,8\}),
\\&(1,3,3),\ (1,3,4),\ (1,3,6),\ (1,3,7),\ (1,3,8),\ (1,3,9).
\endalign$$

{\rm (ii)} The sums
$$\align &p_5+p_5+2p_5,\ p_5+p_5+4p_5,\ p_5+2p_5+2p_5,
\\&p_5+2p_5+4p_5,\ p_5+p_5+5p_5,\ p_5+3p_5+6p_5
\endalign$$
are universal over $\Z$.
\endproclaim

\Remark\ 1.2. As a supplement to Theorem 1.1, the author conjectured in 2009 that $p_5+bp_5+cp_5$ is universal over $\Z$
if the ordered pair $(b,c)$ is among
$$\align&(1,3),\ (1,6),\ (1,8),\ (1,9),\ (1,10),\  (2,3),
\\&(2,6),\ (2,8),\ (3,3),\ (3,4),\ (3,7),\ (3,8),\ (3,9).
\endalign$$
\medskip

There are no universal sums $ap_k+bp_k+cp_k$ over $\N$ with $k>3$. So we are led to study
mixed sums of polygonal numbers.
Though there are infinitely many positive integers which cannot be written as the sum of three squares,
our computation suggests the following somewhat surprising conjecture.

\proclaim{Conjecture 1.3} {\rm (i)} Any $n\in\N$ can be written as the sum of two squares and a pentagonal number,
and as the sum of a triangular number, an even square and a pentagonal number. For $m\in\{3,4,5,\ldots\}$,
if $m\not\eq2,4\pmod8$ then any sufficiently large integer can be written as the sum of two squares and an $m$-gonal number;
also, $m\not\eq2\pmod{32}$ if and only if any sufficiently large integer can be written as the sum of a
triangular number, a square and an $m$-gonal number.

{\rm (ii)} For $n\in\Z^+$ let $r(n)$ denote the number of ways to write $n$ as the sum of
a square, a pentagonal number and a hexagonal number. Then we have $\{r(n):\ n\in\Z^+\}=\Z^+$.
\endproclaim

\Remark\ 1.4. (a) It seems that neither the proof of the Gauss-Legendre theorem nor Cauchy's proof (cf. [MW, N87, N96])
of Fermat's assertion on polygonal numbers
can be adapted to yield a proof of Conjecture 1.3.

(b) We also conjecture that any positive integer $n\not=18$
can be written as the sum of a triangular number, an odd square and a pentagonal number, and as the sum
$p_4(x)+p_5(y)+p_6(-z)$ with $x,y,z\in\N$.

(c) For $n=1,2,3,\ldots$ let $s(n)$ be the least positive integer which can be
written as the sum of a square, a pentagonal number and a hexagonal number in exactly $n$ ways. We guess that
$s(n)/n^2$ has a limit $c\in(0.9,1.2)$ (which is probably $1$) as $n\to+\infty$.

\medskip

Inspired by Cantor's diagonal method in his proof of the uncountability of the set of real numbers, here we pose
a general conjecture on sums of polygonal numbers.

\proclaim{Conjecture 1.5} Let $m\gs 3$ be an integer.
Then every $n\in\N$ can be written in the form
$$p_{m+1}(x_1)+p_{m+2}(x_2)+p_{m+3}(x_3)+r,$$
where $x_1,x_2,x_3\in\N$ and $r\in\{0,\ldots,m-3\}$.
In other words, any $n\in\N$ can be expressed as $p_{m+1}(x_1)+\cdots+p_{2m}(x_m)$ with $x_1,\ldots,x_m\in\N$
and $x_k<2$ for $3<k\ls m$.
\endproclaim

\Remark\ 1.6. Conjecture 1.5 seems challenging and it might remain open for many years in our opinion.
The conjecture in the case $m=3$ is equivalent to the universality of $p_4+p_5+p_6$ (over $\N$).
For $m=3$, $m\in\{4,\ldots,10\}$ and $m\in\{11,\ldots,40\}$ we have verified Conjecture 1.5 for those $n$ not exceeding
$3\times 10^7$, $10^6$, and $10^5$ respectively.
We also guess that for any $m=3,4,\ldots$
all sufficiently large integers have the form $p_m(x)+p_{m+1}(y)+p_{m+2}(z)$ with $x,y,z\in\N$;
for example, our computation via {\tt Mathematica} suggests that
387904 is the largest integer not of the form $p_{20}(x)+p_{21}(y)+p_{22}(z)$.

\medskip

For $m=3,4,\ldots$, clearly
$$8(m-2)p_m(x)+(m-4)^2=((2m-4)x-(m-4))^2.$$
Thus, for given $a,b,c\in\Z^+$ and $i,j,k\in\{3,4,\ldots\}$, that $n=ap_i(x)+bp_j(y)+cp_k(z)$ for some $x,y,z\in\Z$
implies a representation of certain $nq+r$ (with $q$ only depending on $i,j,k$ and $r$ depending on $a,b,c,i,j,k$)
by a ternary quadratic form.

By a deep theory of ternary quadratic forms,
Jones and Pall [JP] proved that for any $n\in\N$ we can write $24n+1$ in the form $x^2+48y^2+144z^2$ with $x,y,z\in\Z$;
equivalently this says that $2p_4+6p_4+p_5$ is universal over $\Z$. They also showed that for $n\in\N$
we can express $24n+5$ as $x^2+y^2+3z^2$ with $x,y,z$ odd and $x\eq y\ (\mo\ 12)$; this is equivalent to the universality
of $p_3+3p_4+2p_5$ over $\Z$. Another result of Jones and Pall [JP] states that for any $n\in\N$ we can write
$24n+29$ in the form $x^2+y^2+3z^2$ with $x,y,z$ odd and $x-y\eq 6\ (\mo\ 12)$;
this is equivalent to the universality of $p_3+6p_3+p_8$ over $\Z$.

Though we are unable to prove Conjecture 1.3, we can show the following result.

\proclaim{Theorem 1.7} {\rm (i)} Let $n\in\N$.
Then $6n+1$ can be expressed in the form $x^2+3y^2+24z^2$ with $x,y,z\in\Z$.
Consequently, there are $x,y,z\in\Z$ such that
$$n=(2x)^2+p_5(y)+p_6(z)=p_4(2x)+p_5(y)+p_3(2z-1),$$
i.e., $n$ can be expressed as the sum of a triangular number, an even square and a generalized pentagonal number.

{\rm (ii)} The sums $p_4+p_4+p_5$ and $p_4+2p_4+p_5$ are universal over $\Z$.

{\rm (iii)}  For any $n\in\N$, we can write $12n+4$ in the form $x^2+3y^2+3z^2$
with $x,y,z\in\Z$ and $2\nmid x$, and write  $12n+4$ and $12n+8$ in the form $3x^2+y^2+z^2$
with $x,y,z\in\Z$ and $2\nmid x$. Also, for any $n\in\N$ the number
 $12n+5$ can be expressed as $x^2+y^2+(6z)^2$ with $x,y,z\in\Z$.
Consequently, $3p_3+3p_4+p_5$ and $p_3+p_4+p_{11}$ are universal over $\Z$.

{\rm (iv)} Let $n\in\N$ and $r\in\{1,9\}$. If $20n+r$ is not a square, then
there are $x,y,z\in\Z$ such that $20n+r=5x^2+5y^2+4z^2$. Consequently, if $20(n+1)+1$ is not a square then
$n=p_3(x)+p_3(y)+p_{12}(z)$ for some $x,y,z\in\Z$.

{\rm (v)}  For any $n\in\N$ with $7n+4$ not squarefree, there are $x,y,z\in\Z$ such that $n=p_3(x)+p_4(y)+p_9(z)$.

{\rm (vi)} If $n\in\N$ is sufficiently large, then we can write it in the form $p_3(x)+p_4(y)+p_{18}(z)$ with $x,y,z\in\Z$,
we can also express $n$ as $p_3(x)+p_4(y)+p_{27}(z)$ with $x,y,z\in\Z$.
Under the GRH (Generalized Riemann Hypothesis for $L$-functions), $p_3+p_4+p_{18}$ and $p_3+p_4+p_{27}$ are universal over $\Z$.
\endproclaim

\Remark\ 1.8. (i) Let $\da\in\{0,1\}$ and $n\in\N$. We conjecture that
if $n\gs\da$ then $n=x^2+y^2+p_5(z)$ for some $x,y,z\in\Z$ with $x\eq \da\ (\mo\ 2)$.
We also conjecture that  for $r\in\{1,9\}$ there are $x,y,z\in\Z$
such that $20n+r=5x^2+5y^2+z^2$ with $z\eq\da\ (\mo\ 2)$ unless
$r=1$, and $\da=n=0$ or ($\da=1\ \&\ n=3)$.

(ii) Our proof of the first assertion in Theorem 1.7(iii) has the same flavor with Fermat's method of infinite descent.

\medskip

As for the determination of those universal sums $p_i+p_j+p_k$ over $\N$, we have the following result.

\proclaim{Theorem 1.9}
Suppose that $p_i+p_j+p_k$ is universal over $\N$ with $3\ls i\ls j\ls k$ and $k\gs 5$. Then
$(i,j,k)$ is among the following $31$ triples:
$$\align&(3,3,5),\ (3,3,6),\ (3,3,7),\ (3,3,8),\ (3,3,10),\ (3,3,12),\ (3,3,17),
\\&(3,4,5),\ (3,4,6),\ (3,4,7),\ (3,4,8),\ (3,4,9),\ (3,4,10),\ (3,4,11),
\\&(3,4,12),\ (3,4,13),\ (3,4,15),\ (3,4,17),\ (3,4,18),\ (3,4,27),
\\&(3,5,5),\ (3,5,6),\ (3,5,7),\ (3,5,8),\ (3,5,9),\ (3,5,11),\ (3,5,13),
\\&(3,7,8),\ (3,7,10),\ (4,4,5),\ (4,5,6).
\endalign$$
\endproclaim

In view of Theorem 1.9, we formulate the following conjecture based on our computation.
\proclaim{Conjecture 1.10} If $(i,j,k)$ is one of the $31$ triples listed in Theorem $1.9$, then
$p_i+p_j+p_k$ is universal over $\N$.
\endproclaim

Besides Theorem 1.9, we also can show that if $3\ls i\ls j\ls k$ and there is a unique nonnegative integer
not in $N(p_i,p_j,p_k)$ then $(i,j,k)$ is among the following 29 triples and the unique number in $\N\sm N(p_i,p_j,p_k)$
does not exceed 468.
$$\align&(3,3,9),\ (3,3,11),\ (3,3,13),\ (3,3,14),\ (3,3,16),\ (3,3,20),\ (3,4,14),
\\& (3,4,19),\ (3,4,20),\ (3,4,22),\ (3,4,29),\ (3,5,10),\ (3,5,12),\ (3,5,14),
\\&(3,5,15),\ (3,5,19),\ (3,5,20),\ (3,5,22),\ (3,5,23),\ (3,5,24),\ (3,5,32),
\\&(3,6,7),\ (3,6,8),\ (3,7,9),\ (3,7,11),\ (3,8,9),\ (3,8,10),\ (4,5,7),\ (4,5,8).
\endalign$$
We also have the following conjecture: If $(i,j,k)$ is among the above 29 triples,
then there is a unique nonnegative integer not in $N(p_i,p_j,p_k)$. In particular,
$$N(p_3,p_4,p_{20})=\N\sm\{468\},\ N(p_3,p_5,p_{32})=\N\sm\{31\}
\ \t{and}\ N(p_4,p_5,p_8)=\N\sm\{19\}.$$

\medskip
Here is another theorem on universal sums over $\N$.

\proclaim{Theorem 1.11} Let $a,b,c\in\Z^+$ with $\max\{a,b,c\}>1$,
and let $i,j,k\in\{3,4,\ldots\}$ with $i\ls j\ls k$ and $\max\{i,j,k\}\gs5$.
Suppose that $(ap_i,bp_j,cp_k)$ is universal $($over $\N)$ with $a\ls b$ if $i=j$, and $b\ls c$ if $j=k$. Then
$(ap_i,bp_j,cp_k)$ is on the following list of $64$ triples:
$$\align&(p_3,p_3,2p_5),\ (p_3,p_3,4p_5),\ (p_3,2p_3,p_5),\ (p_3,2p_3,4p_5),\ (p_3,3p_3,p_5),
\\&(p_3,4p_3,p_5),\ (p_3,4p_3,2p_5),\ (p_3,6p_3,p_5),\ (p_3,9p_3,p_5),\ (2p_3,3p_3,p_5),
\\&(p_3,2p_3,p_6),\ (p_3,2p_3,2p_6),\ (p_3,2p_3,p_7),\ (p_3,2p_3,2p_7),\ (p_3,2p_3,p_8),
\\&(p_3,2p_3,2p_8),\ (p_3,2p_3,p_9),\ (p_3,2p_3,2p_9),\ (p_3,2p_3,p_{10}),\ (p_3,2p_3,p_{12}),
\\&(p_3,2p_3,2p_{12}),\ (p_3,2p_3,p_{15}),\ (p_3,2p_3,p_{16}),\ (p_3,2p_3,p_{17}),\ (p_3,2p_3,p_{23}),
\\&(p_3,p_4,2p_5),\ (p_3,2p_4,p_5),\ (p_3,2p_4,2p_5),\ (p_3,2p_4,4p_5),\ (p_3,3p_4,p_5),
\\&(p_3,4p_4,p_5),\ (p_3,4p_4,2p_5),\ (2p_3,p_4,p_5),\ (2p_3,p_4,2p_5),\ (2p_3,p_4,4p_5),
\\&(2p_3,3p_4,p_5),\ (3p_3,p_4,p_5),\ (p_3,2p_4,p_6),\ (2p_3,p_4,p_6),\ (p_3,p_4,2p_7),
\\&(2p_3,p_4,p_7),\ (p_3,p_4,2p_8),\  (p_3,2p_4,p_8),\ (p_3,3p_4,p_8),\ (2p_3,p_4,p_8),
\\&(2p_3,3p_4,p_8),\ (p_3,p_4,2p_9),\ (p_3,2p_4,p_9),\ (2p_3,p_4,p_{10}),\ (2p_3,p_4,p_{12}),
\\&(p_3,2p_4,p_{17}),\ (2p_3,p_4,p_{17}),\ (p_3,p_5,4p_6),\ (p_3,2p_5,p_6),\ (p_3,p_5,2p_7),
\\&(p_3,p_5,4p_7),\ (p_3,2p_5,p_7),\ (3p_3,p_5,p_7),\ (p_3,p_5,2p_9),\ (2p_3,p_5,p_9),
\\&(p_3,2p_6,p_8),\ (p_3,p_7,2p_7),\ (p_4,2p_4,p_5),\ (2p_4,p_5,p_6).
\endalign$$
\endproclaim

\Remark\ 1.12. It is known (cf. [S07, Lemma 1]) that
$$\{p_3(x)+p_3(y):\ x,y\in\N\}=\{2p_3(x)+p_4(y):\ x,y\in\N\}.\tag1.5$$
So, for $c\in\Z^+$ and $k\in\{3,4,\ldots\}$, we have $N(p_3,p_3,cp_k)=N(2p_3,p_4,cp_k)$
and $Z(p_3,p_3,cp_k)=Z(2p_3,p_4,cp_k)$.

\medskip

Here we pose the following conjecture based on our computation.

\proclaim{Conjecture 1.13} All the $64$ triples listed in Theorem $1.11$
are universal over $\N$.
\endproclaim

We are unable to show the universality over $\N$ for any of the $31+64=95$ triples in Theorems 1.9 and 1.11,
but we can prove that many of them
are universal over $\Z$.
In light of (1.3) and (1.5), when we consider the universality of $ap_i+bp_j+cp_k$ over $\Z$,
we may ignore those triples $(ap_i,bp_j,cp_k)$ with
$(ap_i,bp_j)=(2p_3,p_4)$ or $6\in\{i,j,k\}$. This reduces the 31+64 triples to $(31-4)+(64-16)=75$ essential triples.
By Theorem 1.7, the sums
$$p_3+p_4+p_5,\ p_3+4p_4+p_5,\ p_4+p_4+p_5,\ p_4+2p_4+p_5,\ p_3+p_4+p_{11}$$
are universal over $\Z$. So there are $75-5=70$ essential triples left.

A positive ternary quadratic form $Q(x,y,z)=ax^2+by^2+cz^2+dyz+exz+fxy$
with $a,b,c,d,e,f\in\Z$ is said to be {\it regular}
if it represents an integer $n$ (i.e., $Q(x,y,z)=n$ for some $x,y,z\in\Z$) if and only if
it locally represents $n$ (i.e., for any prime $p$ the equation $Q(x,y,z)=n$
has integral solutions in the $p$-adic field $\Q_p$;
in other words, for any $m\in\Z^+$ the congruence $Q(x,y,z)\eq n\ (\mo\ m)$ is solvable over $\Z$).
A full list of positive regular ternary quadratic forms was given in [JKS].
There are totally 102 regular forms $ax^2+by^2+cz^2$ with $1\ls a\ls b\ls c$ and $\gcd(a,b,c)=1$; for each of them
those positive integers not represented by the form were described explicitly in Dickson [D39, pp.112-113].

By applying the theory of ternary quadratic forms and using some lemmas of ourselves,
we are able to deduce the following result.

\proclaim{Theorem 1.14} $ap_i+bp_j+cp_k$ is universal over $\Z$ if $(ap_i,bp_j,cp_k)$ is one of the following $35$ essential triples:
$$\align&(p_3,p_3,p_5),\ (p_3,p_3,2p_5),\ (p_3,p_3,4p_5),\ (p_3,2p_3,p_5),\ (p_3,2p_3,4p_5),
\\&(p_3,3p_3,p_5),\ (p_3,4p_3,p_5),\ (p_3,4p_3,2p_5),\ (p_3,6p_3,p_5),\ (2p_3,3p_3,p_5),
\\&(p_3,p_4,2p_5),\ (p_3,2p_4,p_5),\ (p_3,3p_4,p_5),\ (p_3,2p_4,2p_5),\ (2p_3,3p_4,p_5),
\\&(3p_3,p_4,p_5),\ (p_3,p_5,p_5),\ (p_3,p_3,p_7),\ (p_3,2p_3,2p_7),\ (p_3,p_4,p_7),
\\&(p_3,p_5,2p_7),\ (p_3,p_3,p_8),\ (p_3,2p_3,p_8),\ (p_3,p_4,p_8),\ (p_3,p_4,2p_8),
\\&(p_3,2p_4,p_8),\ (p_3,3p_4,p_8),\ (2p_3,3p_4,p_8),\ (p_3,p_5,p_8),\ (p_3,p_5,p_9),
\\&(p_3,p_3,p_{10}),\ (p_3,2p_3,p_{10}),\ (p_3,p_4,p_{10}),\ (p_3,p_7,p_{10}),\ (p_3,p_4,p_{12}).
\endalign$$
\endproclaim

For the following $70-35=35$ remaining essential triples $(ap_i,bp_j,cp_k)$, we have not yet proved
the universality of $ap_i+bp_j+cp_k$ over $\Z$.
$$\align&(p_3,9p_3,p_5),\ (p_3,2p_4,4p_5),\ (p_3,4p_4,2p_5),\ (p_3,2p_3,p_7),\ (p_3,p_4,2p_7),
\\&(p_3,p_5,p_7),\ (p_3,p_5,4p_7),\ (p_3,2p_5,p_7),\ (3p_3,p_5,p_7),\ (p_3,p_7,2p_7),
\\&(p_3,2p_3,2p_8),\ (p_3,p_7,p_8),\ (p_3,2p_3,p_9),\ (p_3,2p_3,2p_9),\ (p_3,p_4,p_9),
\\&(p_3,p_4,2p_9),\ (p_3,2p_4,p_9),\ (p_3,p_5,2p_9),\ (2p_3,p_5,p_9),\ (p_3,p_5,p_{11}),
\\&(p_3,p_3,p_{12}),\ (p_3,2p_3,p_{12}),\ (p_3,2p_3,2p_{12}),\ (p_3,p_4,p_{13}),\ (p_3,p_5,p_{13}),
\\&(p_3,2p_3,p_{15}),\ (p_3,p_4,p_{15}),\ (p_3,2p_3,p_{16}),\ (p_3,p_3,p_{17}),\ (p_3,2p_3,p_{17}),
\\&(p_3,p_4,p_{17}),\ (p_3,2p_4,p_{17}),\ (p_3,p_4,p_{18}),\ (p_3,2p_3,p_{23}),\ (p_3,p_4,p_{27}).
\endalign$$

\medskip

It is easy to see that
$$\align &n=p_3(x)+p_4(y)+p_{17}(z)
\\\iff&120n+184=15(2x+1)^2+120y^2+(30z-13)^2
\endalign$$
and
$$\align &n=p_3(x)+2p_3(y)+p_{23}(z)
\\\iff&168n+424=21(2x+1)^2+42(2y+1)^2+(42z-19)^2.
\endalign$$
However, when $w^2\eq 13^2\ (\mo\ 120)$ we may not have $w\eq\pm 13\ (\mo\ 30)$ since $23^2\eq 13^2\ (\mo\ 120)$;
similarly, when $w^2\eq 19^2\ (\mo\ 168)$ we may not have $w\eq\pm 19\ (\mo\ 42)$ since $19^2\eq 5^2\ (\mo\ 168)$.
Thus it seems quite difficult to show that the sums
$$p_3+p_3+p_{17},\, p_3+2p_3+p_{17},\,p_3+p_4+p_{17},\,p_3+2p_4+p_{17},\,p_3+2p_3+p_{23}$$
are universal over $\Z$, let alone their universality over $\N$.

The study of some of the 35 remaining essential triples leads us to pose the following conjecture.

\proclaim{Conjecture 1.15} If $a\in\Z^+$ is not a square, then sufficiently large integers relatively prime to $a$
can be written in the form
$p+ax^2$ with $p$ prime and $x\in\Z$, i.e., the set $S(a)$ given by
$$\{n>1:\, \gcd(a,n)=1,\ \t{and}\ n\not=p+ax^2\ \t{for any prime}\ p \ \t{and}\ x\in\Z\}\tag1.6$$
is finite.
In particular,
$$\gather S(6)=\em,\ \ S(12)=\{133\},\ \ S(30)=\{121\},
\\ S(3)=\{4,\,28,\,52,\,133,\,292,\,892,\,1588\},\, S(18)=\{187,\,1003,\,5777,\,5993\},
\\ S(24)=\{25,\,49,\,145,\,385,\,745,\,1081,\,1139,\,1561,\,2119,\,2449,\,5299\}.
\endgather$$
Also, we have the precise values of $M(a)=\max S(a)$ for some other values of $a$:
$$M(5)=270086,\ M(7)=150457,\ M(8)=39167,\ M(10)=18031,\ M(11)=739676.$$
\endproclaim
\Remark\ 1.16. According to [D99a, p.\,424], in 1752 Goldbach
asked whether any odd integer $n>1$ has the form $p+2x^2$,
and $5777,\,5993\in S(2)$ was found by M. A. Stern and his students in 1856. It seems that $S(2)=\{5777,\, 5993\}$.
See [S09] for the author's conjectures on sums of primes and triangular numbers.

\proclaim{Theorem 1.17} Under Conjecture $1.15$,
the sums $p_3+2p_4+p_9$ and $p_3+p_4+p_{13}$ are  universal over $\Z$.
\endproclaim

\medskip

Our following conjecture, together with Theorems 1.9 and 1.11, Conjectures 1.10 and 1.13, and (1.3),
describes all universal sums $p_i+p_j+p_k$ over $\Z$.

\proclaim{Conjecture 1.18} For $i,j,k\in\{3,4,\ldots\}\sm\{6\}$ with $i\ls j\ls k$ and $k\gs5$,
$p_i+p_j+p_k$ is universal over $\Z$ but not universal over $\N$, if and only if $(i,j,k)$ is among the following list:
$$\align&(3,3,k)\ (k\in\{9,11,13,14,15,16,20,23,24,25,26,29,32,33,35\}),
\\& (3,4,k)\ (k\in\{14,16,19,20,21,22,23,24,26,29,30,32,33,35,37\}),
\\&(3,5,k)\ (k\in[10,68]\sm\{11,13,26,34,36,44,48,56,59,60,64\}),
\\&(3,7,k)\ (k\in[7,54]\sm\{8,10,42,51\}),
\ (3,8,k)\ (k\in[8,16]),
\\&(3,9,k)\ (k\in[10,17]\sm\{13\}),\ (3,10,k)\, (k=11,\ldots,22),
\\&(3,11,k)\ (k\in[12,23]\sm\{18,19\}),\ (3,12,k)\ (k\in[13,27]\sm\{19,21,22,25\});
\endalign$$
$$\align&(4,4,7), (4,4,8), (4,4,10),(4,5,5), (4,5,k)\,(k\in[7,37]),
\\&(4,7,k)\ (k\in[7,20]\cup\{23,25,26,27,29,38,41,44\}),
\\&(4,8,k)\ (k\in[8,22]\sm\{13,14,20\}),\ (4,9,k)\ (k=9,\ldots,15),
\\&(4,10,11),(4,10,12),(4,10,14);
\endalign$$
$$\align&(5,5,k)\ (k\in[5,71]\sm\{6,42,45,50,56,58,59,61,64,67,69,70\}),
\\&(5,7,k)\, (k\in[7,41]\sm\{17\}),\ (5,8,k)\, (k\in[8,28]\sm\{14,21,25\}),
\\&(5,9,k)\, (k\in[9,41]\sm\{26\}),\ (5,10,k)\, (k=10,\ldots,24),
\\&(5,11,k)\, (k\in[11,32]\sm\{15,18,25,26,29\}),
\\&(5,12,k)\, (k\in[12,76]\sm\{30,35,36,42,49,52,53,55,56,64,71,73,74\}),
\\& (5,13,k)\,(k\in[13,33]\sm\{30\});
\endalign$$
$$\align
&(7,7,8),\ (7,7,9),\ (7,7,11),\ (7,8,k)\, (k\in[8,14]),\ (7,9,k)\, (k\in[9,15]),
\\&(7,10,k)\, (k\in[11,19]\sm\{13,18\}),
\ (7,11,k)\, (k\in[11,19]\sm\{14,18\}),
\\&(7,12,k)\, (k\in[13,23]\sm\{18,21,22\}),\ (7,13,14),\ (7,13,15),\ (7,13,17).
\endalign$$
\endproclaim
\Remark\ 1.19. We also have a conjecture describing all those universal sums $ap_i+bp_j+cp_k$ over $\Z$
 with $a,b,c\in\Z^+$.
There are only finitely many such sums (including $18p_3+p_4+p_5$, $p_4+p_5+20p_5$ and $p_8+3p_8+p_{10}$).

\medskip

Via the theory of quadratic forms, we can prove that
many of those $p_i+p_j+p_k$ with $(i,j,k)$ listed in Conjecture 1.18 are indeed universal over $\Z$.
Here we include few interesting cases in the following theorem.

\proclaim{Theorem 1.20} The sums
$$p_3+p_7+p_7,\ p_4+p_5+p_5,\ \ p_4+p_4+p_8,\ p_4+p_8+p_8,\ \ p_4+p_4+p_{10},\ p_5+p_5+p_{10}$$
are universal over $\Z$ but not universal over $\N$.
\endproclaim

\Remark\ 1.21. Let $\da\in\{0,1\}$ and $n\in\{\da,\da+1,\da+2,\ldots\}$.
We conjecture that $n=x^2+p_5(y)+p_5(z)$ for some $x,y,z\in\Z$ with $x\eq\da\ (\mo\ 2)$.
\medskip

Among those triples $(i,j,k)$ listed in Conjecture 1.18, the one with $i+j+k$ maximal is
$(5,12,76)$. By Conjecture 1.18 the sum $p_5+p_{12}+p_{76}$ should be universal over $\Z$. Observe that
$n=p_5(x)+p_{12}(y)+p_{76}(z)$ if and only if
$$4440(n+9)+2657=185(6x-1)^2+888(5y-2)^2+120(37z-18)^2.$$
It seems quite difficult to prove that for any $n=9,10,\ldots$ the equation
$$4440n+2657=185x^2+888y^2+120z^2\tag1.7$$
has integral solutions.

In Sections 2 and 3 we shall show Theorems 1.1 and 1.7 respectively.
Theorem 1.14 will be proved in Section 4. Section 5 is devoted to the proofs of Theorems 1.17 and 1.20.
We will give two auxiliary results in Section 6 and show Theorems 1.9 and 1.11 in Sections 7-8 respectively.
To make this paper easily understood for general readers, we try our best to avoid using specific tools only known to experts at quadratic forms.

\heading{2. Proof of Theorem 1.1}\endheading

\medskip
\noindent{\it Proof of Theorem {\rm 1.1(i)}}. As $k\gs 4$, $p_k(3)=3k-3\gs9$ and $p_k(-3)=6k-15\gs 9$, only the following numbers
$$p_k(0)=0,\ p_k(1)=1,\ p_k(-1)=k-3,\ p_k(2)=k,\ p_k(-2)=3k-8$$
can be elements of $p_k(\Z)$ smaller than 8. By  $1\in Z(ap_k,bp_k,cp_k)$ (and $a\ls b\ls c$),
we get $a=1$. Since $2,3\in Z(p_k,bp_k,cp_k)$, if $b>2$ then $k=5$ and $b=3$.

\medskip

{\it Case} 1. $b=1$.

 Observe that
$$p_5(\Z)=\bg\{\f{3n^2\pm n}2:\
n=0,1,\ldots\bg\}=\{0,1,2,5,7,12,15,22,26,35,40,\ldots\}$$
and thus $p_5(\Z)+p_5(\Z)$ does not contain 11. If $k=5$, then $c$ cannot be greater than 11. It is easy
to verify that $25\not\in Z(p_5,p_5,7p_5)$ and $43\not\in Z(p_5,p_5,11p_5)$.
So $c\in[1,10]\sm\{7\}$ when $k=5$.

 Now assume that $k\not=5$. Then
$3\not\in p_k(\Z)+p_k(\Z)$. By $3\in Z(p_k,p_k,cp_k)$, we must have $c\ls 3$.
Observe that $7\not\in Z(p_4,p_4,p_4)$, $14\not\in Z(p_4,p_4,2p_4)$
and $6\not\in Z(p_4,p_4,3p_4)$. So $k\gs 7$. It is easy to verify that
$$Z(p_k,p_k,p_k)\cap[4,7]=\{k-3,\,k-2,\,k-1,\,k\}\cap[4,7].$$
Thus, if $c=1$ then $k=7$. But $10\not\in Z(p_7,p_7,p_7)$, so $c\in\{2,3\}$.
Observe that
$$Z(p_k,p_k,2p_k)\cap[5,7]=\{k-3,\,k-2,\,k-1,\,k\}\cap[5,7].$$
If $c=2$, then $k\in\{7,8\}$. But $23\not\in Z(p_7,p_7,2p_7)$ and
$14\not\in Z(p_8,p_8,2p_8)$, therefore $c=3$. Clearly
$$Z(p_k,p_k,3p_k)\cap[6,7]=\{k-3,\,k-2,\,k\}\cap[6,7].$$
So $k=9$. As $8\not\in Z(p_9,p_9,3p_9)$, we get a contradiction.

\medskip

{\it Case} 2. $b\in\{2,3\}$.

For $b=2,3$ it is easy to see that $b+6\not\in p_5(\Z)+bp_5(\Z)$.
If $k=5$, then  $b+6\in Z(p_5,bp_5,cp_5)$ and hence $c\ls b+6$.
Observe that $18\not\in Z(p_5,2p_5,5p_5)$, $27\not\in Z(p_5,2p_5,7p_5)$
and $19\not\in Z(p_5,3p_5,5p_5)$.

 Now suppose that $k\not=5$. Then we must have $b=2$.
 If $k=4$, then by
 $5\in Z(p_4,2p_4,cp_4)$ we get $c\ls 5$. But, for $c=2,3,4,5$
 the set $Z(p_4,2p_4,cp_4)$ does not contain 7, 10, 14, 10 respectively. So $k\gs 7$.
Note that $c(k-3)\gs2(k-3)>7$ and
 $$Z(p_k,2p_k,cp_k)\cap[0,7]\subseteq\{0,\,1,\,2,\,3,\,k-3,\,k-1,\,k\}+\{0,c\}.$$
By $4\in Z(p_k,2p_k,cp_k)$, we have $k=7$ or $c\ls 4$.
By $5\in Z(p_k,2p_k,cp_k)$, we have $k=8$ or $c\ls 5$.
Therefore $c\ls 4$, or $c=5$ and $k=7$. For $c=2,3,4,5$ the set $Z(p_7,2p_7,cp_7)$
does not contain 19, 31, 131, 10 respectively. So $c\in\{2,3,4\}$ and $k\gs8$. Since
$$6,7\in\{0,\,1,\,2,\,3,\,k-3,\,k-1\}+\{0,\,c\},$$
if $c\in\{2,3\}$ then $c=3$ and $k\in\{8,10\}$. But $9\not\in Z(p_8,2p_8,3p_8)$
and $8\not\in Z(p_{10},2p_{10},3p_{10})$, so we must have $c=4$.
For $k=8, 9,10,11,12,\ldots$ we have $a_k\not\in Z(p_k,2p_k,4p_k)$, where
$$a_8=13,\ a_9=14,\ a_{11}=9,\ a_{10}=a_{12}=a_{13}=\cdots=8.$$ So we get a contradiction.

 In view of the above we have proved part (i) of Theorem 1.1. \qed

 \proclaim{Lemma 2.1} Let $w=x^2+my^2$ be a positive integer
with $m\in\{2,5,8\}$ and $x,y\in\Z$. Then we can write $w$ in the form $u^2+mv^2$ with $u,v\in\Z$
such that $u$ or $v$ is not divisible by $3$.
\endproclaim
\Proof. Write $x=3^kx_0$ and $y=3^ky_0$ with $k\in\N$,
where $x_0$ and $y_0$ are integers not all divisible by 3.
Then $w=9^k(x_0^2+my_0^2)$.

Suppose that $a$ and $b$ are integers not all divisible by 3. It is easy to check the following identities:
$$\align9(a^2+2b^2)=&(1^2+2\times2^2)(a^2+2b^2)=(a\pm 4b)^2+2(2a\mp b)^2,
\\9(a^2+5b^2)=&(2^2+5\times1^2)(a^2+5b^2)=(2a\pm 5b)^2+5(a\mp 2b)^2,
\\9(a^2+8b^2)=&(1^2+8\times1^2)(a^2+8b^2)=(a\pm 8b)^2+8(a\mp b)^2.
\endalign$$
For each $k=1,2,4,8$ we cannot have $a+kb\eq a-kb\eq0\ (\mo\ 3)$. Thus
we can rewrite $9(a^2+mb^2)$ in the form $c^2+md^2$ with $c$ or $d$ not divisible by 3.

Applying the above process repeatedly, we finally get that $w=u^2+mv^2$ for some $u,v\in\Z$ not all divisible by 3.
This concludes the proof. \qed

\medskip
\Remark\ 2.2. Lemma 2.1 in the case $m=2$ first appeared in the middle of a proof given on page 173 of [JP].

\medskip
\noindent{\it Proof of Theorem {\rm 1.1(ii)}}. Let $n$ be any nonnegative integer.
If $w\in\Z$ is relatively prime to 6, then $w$ or $-w$ has the form $6x-1$ with $x\in\Z$.
Thus, given $b,c\in\Z^+$ we have
$$\align &n=p_5(x)+bp_5(z)+cp_5(z)\ \t{for some}\ x,y,z\in\Z
\\\iff& \bar n=(6x-1)^2+b(6y-1)^2+c(6z-1)^2\ \t{for some}\ x,y,z\in\Z
\\\iff&\bar n=x^2+by^2+cz^2\
 \t{for some integers}\ x,y,z\ \t{relatively prime to}\ 6,
\endalign$$
where $\bar n=24n+b+c+1$. Below we will often appeal to this basic fact.

 (a) By the Gauss-Legendre theorem, there are $u,v,w\in\Z$ such that
 $12n+2=u^2+v^2+w^2$. As $u^2+v^2+w^2\eq2\ (\mo\ 3)$, exactly one of $u,v,w$ is divisible by $3$.
 Without loss of generality, we assume that $3\mid u$ and $3\nmid vw$.
 Clearly $u,v,w$ cannot have the same parity.
 Without loss of generality, we suppose that $v\not\eq u\ (\mo\ 2)$.
 Note that both $u\pm v$ and $w$ are relatively prime to $6$. Since $24n+4=(u+v)^2+(u-v)^2+2w^2$,
 there are $x,y,z\in\Z$ such that $n=p_5(x)+p_5(y)+2p_5(z)$.

(b) Write $24n+6=9^kn_0$ with $k,n_0\in\N$ and $9\nmid n_0$. Obviously $n_0\eq 6\ (\mo\ 8)$.
By the Gauss-Legendre theorem, $n_0=x_0^2+y_0^2+z_0^2$ for some $x_0,y_0,z_0\in\Z$. Clearly
$x_0,y_0,z_0$ are not all divisible by 3. If $x,y,z\in\Z$ and $3\nmid x$, then $x'=\ve x\not\eq 2y+2z\ (\mo\ 3)$
for a suitable choice of $\ve\in\{1,-1\}$, hence by the identity
$$(2y+2z-x')^2+(2x'+2z-y)^2+(2x'+2y-z)^2=9((x')^2+y^2+z^2)$$
(cf. (1.4)) we see that $9(x^2+y^2+z^2)$ can be written as $u^2+v^2+w^2$ with $3\nmid u$.
Thus, there are $u,v,w\in\Z$ such that $24n+6=u^2+v^2+w^2$ and not all of $u,v,w$ are multiples of 3.
Since $u^2+v^2+w^2\eq0\ (\mo\ 3)$, none of $u,v,w$ is divisible by $3$. Without loss of generality, we assume that
$w$ is even. As $u^2+v^2\eq 6\ (\mo\ 4)$, we get $u\eq v\eq1\ (\mo\ 2)$. Note that
$w^2\eq 6-u^2-v^2\eq 4\ (\mo\ 8)$ and hence $4\nmid w$. So $u,v,w/2$ are all relatively prime to 6.
Thus, by $24n+6=u^2+v^2+4(w/2)^2$, we have $n=p_5(x)+p_5(y)+4p_5(z)$ for some $x,y,z\in\Z$.

(c) By the Gauss-Legendre theorem, we can write $24n+5=u^2+s^2+t^2$ with $u,s,t\in\Z$ and $s\eq t\ (\mo\ 2)$.
Set $v=(s+t)/2$ and $w=(s-t)/2$. Then $24n+5=u^2+(v+w)^2+(v-w)^2=u^2+2v^2+2w^2$.
As $u^2\not\eq 5\ (\mo\ 3)$, without loss of generality we may assume that $3\nmid w$.
Note that $2\nmid u$.
By Lemma 2.1, $u^2+2v^2=a^2+2b^2$ for some $a,b\in\Z$ with $a,b\in\Z$ not all divisible by 3.
Since $a^2+2b^2=u^2+2v^2\eq 5-2w^2\eq0\ (\mo\ 3)$, both $a$ and $b$ are relatively prime to 3.
As $24n+5=a^2+2b^2+2w^2$, we have $a\eq1\ (\mo\ 2)$ and $2(b^2+w^2)\eq5-a^2\eq 4\ (\mo\ 8)$.
Therefore $a,b,w$ are all relatively prime to 6. So
$n=p_5(x)+2p_5(y)+2p_5(z)$ for some $x,y,z\in\Z$.

(d) By the Gauss-Legendre theorem, $48n+14=s^2+t^2+(2w)^2$ for some $s,t,w\in\Z$ with $s\eq t\ (\mo\ 2)$.
Set $u=(s+t)/2$ and $v=(s-t)/2$. Then $48n+14=(u+v)^2+(u-v)^2+4w^2$ and hence
$24n+7=u^2+v^2+2w^2$. As $2w^2\not\eq 7\ (\mo\ 3)$, without loss of generality we may assume that $3\nmid u$.
Clearly $v^2+2w^2>0$ since $u^2\not\eq 7\ (\mo\ 8)$. By Lemma 2.1, $v^2+2w^2=a^2+2b^2$ for some $a,b\in\Z$ with
$a$ or $b$ not divisible by 3. Note that $a^2+2b^2\eq 7-u^2\eq0\ (\mo\ 3)$ and hence both $a$ and $b$ are relatively
prime to 3. As $24n+7=u^2+a^2+2b^2$, we have $u\not\eq a\ (\mo\ 2)$ and hence $u^2+a^2\eq 1\ (\mo\ 4)$.
Thus $2b^2\eq 7-1\eq 2\ (\mo\ 4)$ and hence $2\nmid b$. Since $u^2+a^2\eq7-2b^2\eq 5\ (\mo\ 8)$,
one of $u$ and $a$ is odd and the other is congruent to 2 mod 4. So, there are $x,y,z\in\Z$ such that
$$24n+7=u^2+a^2+2b^2=(6x-1)^2+2(6y-1)^2+(2(6z-1))^2$$
and hence $n=p_5(x)+2p_5(y)+4p_5(z)$.
\medskip

(e) Recall that $p_3+p_3+5p_3$ is universal as obtained by Liouville.
So there are $r,s,t\in\Z$ such that $3n=p_3(r)+p_3(s)+5p_3(t)$ and hence
$24n+7=u^2+v^2+5w^2$ where $u=2r+1$, $v=2s+1$ and $w=2t+1$. As $5w^2\not\eq7\ (\mo\ 3)$,
without loss of generality we may assume that $3\nmid u$. As $u^2\not\eq 7\ (\mo\ 8)$, we have $v^2+5w^2>0$.
In light of Lemma 2.1, $v^2+5w^2=a^2+5b^2$ for some $a,b\in\Z$ with
$a$ or $b$ not divisible by 3. Since $a^2+5b^2\eq 7-u^2\eq0\ (\mo\ 3)$, both $a$ and $b$ are relatively
prime to 3. As $a^2+5b^2=v^2+5w^2\eq7-u^2\eq6\ (\mo\ 8)$, we get $a\eq b\eq1\ (\mo\ 2)$.
Note that $a,b,u$ are all relatively prime to 6. So, by $24n+7=u^2+a^2+5b^2$, we have
$n=p_5(x)+p_5(y)+5p_5(z)$ for some $x,y,z\in\Z$.
\medskip

(f) By Theorem 1.7(iii), $12n+5=u^2+v^2+(6w)^2$ for some $u,v,w\in\Z$.
(This will be proved in the next section without appeal to Theorem 1.1.)
It follows that $24n+10=(u+v)^2+(u-v)^2+2(6w)^2$. As $u\not\eq v\ (\mo\ 2)$, both $u+v$ and $u-v$
are odd. Since $(u+v)^2+(u-v)^2\eq 10\eq1\ (\mo\ 3)$, we can write $\{u+v,u-v\}$ as $\{x,3r\}$
with $\gcd(x,6)=1$ and $r\eq1\ (\mo\ 2)$. Thus $24n+10=x^2+9(r^2+8w^2)$. Clearly $x^2\not\eq 10\ (\mo\ 8)$ and hence
$r^2+8w^2>0$. By Lemma 2.1 with $m=8$, there are integers $s$ and $t$ not all divisible by 3 such that
$r^2+8w^2=s^2+8t^2$. Since we cannot have $s+2t\eq s-2t\eq0\ (\mo\ 3)$,
without loss of generality we assume that $s-2t\not\eq0\ (\mo\ 3)$.
(If $s+2t\not\eq0\ (\mo\ 3)$ then we use $-t$ instead of $t$.)
Set $y=s+4t$ and $z=s-2t$. Then $y\eq z\not\eq0\ (\mo\ 3)$.
Note also that $y\eq z\eq s\eq r\eq1\ (\mo\ 2)$. Thus $\gcd(xyz,6)=1$.
Observe that
$$\align 24n+10=&x^2+9(s^2+8t^2)=x^2+(3s)^2+2(6t)^2
\\=&x^2+(y+2z)^2+2(y-z)^2=x^2+3y^2+6z^2.
\endalign$$
So $n=p_5(a)+3p_5(b)+6p_5(c)$ for some $a,b,c\in\Z$.

\medskip

By the above we have completed the proof of Theorem 1.1(ii). \qed

 \heading{3. Proof of Theorem 1.7}\endheading

\proclaim{Lemma 3.1} Let $w\in\N$.
Then $w$ can be written in the form $3x^2+6y^2$ with $x,y\in\Z$, if and only if
$3\mid w$ and $w=u^2+2v^2$ for some $u,v\in\Z$.
\endproclaim
\Proof. When $w=3x^2+6y^2=3(x^2+2y^2)$, we clearly have $w=(x+2y)^2+2(x-y)^2$.

Now suppose that there are $u,v\in\Z$ such that $u^2+2v^2=w\eq0\ (\mo\ 3)$.
Since $u^2\eq v^2\ (\mo\ 3)$, without loss of generality we assume that $u\eq v\ (\mo\ 3)$.
Set $x=(u+2v)/3$ and $y=(u-v)/3$. Then
$$w=u^2+2v^2=(x+2y)^2+2(x-y)^2=3x^2+6y^2.$$

By the above the desired result follows. \qed

\proclaim{Lemma 3.2} Let $n$ be any nonnegative integer. Then
we have
$$\aligned&|\{(x,y)\in\Z^2:\ x^2+3y^2=8n+4\ \t{and}\ 2\nmid x\}|
\\&\ \ =\f23|\{(x,y)\in\Z^2:\ x^2+3y^2=8n+4\}|.
\endaligned\tag3.1$$
\endproclaim
\Proof. Clearly (3.1) is equivalent to $|S_1|=2|S_0|$, where
$$S_r=\{(x,y)\in\Z^2:\ x^2+3y^2=8n+4\ \t{and}\ x\eq y\eq r\ (\mo\ 2)\}.$$
If $(2x,2y)\in S_0$, then $2n+1=x^2+3y^2$ and $x\not\eq y\ (\mo\ 2)$, hence
$(x+3y,\pm(x-y))\in S_1$ since
$$(x+3y)^2+3(x-y)^2=4(x^2+3y^2)=4(2n+1)=8n+4.$$
On the other hand, if $(u,v)\in S_1$ with $u\eq \pm v\ (\mo\ 4)$, then
we have
$$u=x+3y, \ v=\pm(x-y),\ \t{and}\ (2x,2y)\in S_0,$$
where $x=(u\pm3v)/4$ and $y=(u\mp v)/4$.
For $(2x,2y),(2s,2t)\in S_0$, if
$$\{(x+3y,x-y),(x+3y,y-x)\}=\{(s+3t,s-t),(s+3t,t-s)\}$$
then $s-t\eq s+3t=x+3y\eq x-y\not\eq y-x\ (\mo\ 4)$, hence $s-t=x-y$ and $(s,t)=(x,y)$.
Therefore, we do have $|S_1|=2|S_0|$ as desired. \qed

\proclaim{Lemma 3.3} Let $n\in\N$ with $6n+1$ not a square. Then, for any $\da\in\{0,1\}$,
we can write $6n+1$ in the form $x^2+3y^2+6z^2$ with $x,y,z\in\Z$ and $x\eq \da\ (\mo\ 2)$.
\endproclaim
\Proof. By the Gauss-Legendre theorem, there are $r,s,t\in\Z$ such that $12n+2=(2r)^2+s^2+t^2$
and hence $6n+1=2r^2+u^2+v^2$ where $u=(s+t)/2$ and $v=(s-t)/2$.  Since $2r^2\not\eq1\ (\mo\ 3)$,
either $u$ or $v$ is not divisible by 3. Without loss of generality we assume that $3\nmid v$.
As $2r^2+u^2\eq0\ (\mo\ 3)$, by Lemma 3.1 we have $2r^2+u^2=3a^2+6b^2$ for some $a,b\in\Z$.
So, if $v\eq\da\ (\mo\ 2)$, then the desired result follows.

Now suppose that $v\not\eq\da\ (\mo\ 2)$. As $6n+1$ is not a square, $2r^2+u^2>0$. By Lemma 2.1, we can write $u^2+2r^2$
in the form $c^2+2d^2$ with $c$ or $d$ not divisible by 3.
Since $3\mid (u^2+2r^2)$, both $c$ and $d$ are relatively prime to 3.
Note that $6n+1=c^2+2d^2+v^2$ and $c\eq \da\ (\mo\ 2)$ (since $c\not\eq v\ (\mo\ 2)$).
As $2d^2+v^2\eq 1-c^2\eq0\ (\mo\ 3)$, by Lemma 3.1 we can write $2d^2+v^2$ in the form $3y^2+6z^2$ with $y,z\in\Z$.
So the desired result follows. \qed
\medskip

\Remark\ 3.4. We conjecture that for any positive integer $n$ with $6n+1$ a square we also have $6n+1=x^2+3y^2+6z^2$ for some $x,y,z\in\Z$ with $x$ even.
can be replaced by $n\gs1-\da$.

\proclaim{Lemma 3.5} Let $w=x^2+4y^2$ be a positive integer
with $x,y\in\Z$. Then we can write $w$ in the form $u^2+4v^2$ with $u,v\in\Z$
 not all divisible by $5$.
\endproclaim
\Proof. Write $x=5^kx_0$ and $y=5^ky_0$ with $k\in\N$,
where $x_0$ and $y_0$ are integers not all divisible by 5.
Then $w=5^{2k}(x_0^2+4y_0^2)$. If $a$ and $b$ are integers not all divisible by 5, then
we cannot have $a+4b\eq a-4b\eq0\ (\mo\ 5)$, hence by the identity
$$5(a^2+4b^2)=(1^2+4\times1^2)(a^2+5b^2)=(a\pm 4b)^2+4(a\mp b)^2$$
we can rewrite $5(a^2+4b^2)$ in the form $c^2+4d^2$ with $c$ or $d$ not divisible by 5.
Applying the process repeatedly, we finally get that $w=u^2+4v^2$ for some $u,v\in\Z$ not all divisible by 5.
This ends the proof. \qed

\proclaim{Lemma 3.6} Let $w=x^2+7y^2$ be a positive integer
with $x,y\in\Z$ and $8\mid w$. Then we can write $w$ in the form $u^2+7v^2$ with $u$ and $v$ both odd.
\endproclaim
\Proof. Write $x=2^kx_0$ and $y=2^ky_0$ with $k\in\N$,
where $x_0$ and $y_0$ are integers not all even.
Then $w=4^{k}(x_0^2+7y_0^2)$. If $2\mid x_0y_0$, then
$k\gs 2$ since $2\nmid (x_0^2+7y_0^2)$, hence
$$w=4^{k-2}\l((3x_0+7y_0)^2+7(x_0-3y_0)^2\r)$$
with $3x_0+7y_0$ and $x_0-3y_0$ odd. Thus it remains to prove that
if $s$ and $t$ are odd then $4(s^2+7t^2)=u^2+7v^2$ for some odd integers $u$ and $v$.

Let $s$ and $t$ be odd integers. Then $s$ is congruent to $t$ or $-t$ mod 4.
Without loss of generality we assume that $s\eq t\ (\mo\ 4)$. Then
$$4(s^2+7t^2)=\l(\f{3s+7t}2\r)^2+7\l(\f{s-3t}2\r)^2,$$
where $(3s+7t)/2=3(s+t)/2+2t$ and $(s-3t)/2=(s+t)/2-2t$ are odd integers.

Combining the above we obtain the desired result. \qed

\medskip
\noindent{\it Proof of Theorem 1.7}. (i) We first show that
$6n+1=x^2+3y^2+24z^2$ for some $x,y,z\in\Z$.
If $6n+1=m^2$ for some $m\in\Z$, then $6n+1=m^2+3\times 0^2+24\times0^2$.
Now assume that $6n+1$ is not a square. By Lemma 3.3, there are $x,y,z\in\Z$ with $x\not\eq n\ (\mo\ 2)$
such that $6n+1=x^2+3y^2+6z^2$. Observe that $y\eq x+1\eq n\ (\mo\ 2)$ and
$x^2+3y^2\eq(n+1)^2+3n^2\eq 6n+1\ (\mo\ 4)$. So $z$ must be even.

 By the above, there are $x,y,z\in\Z$ such that $6n+1=x^2+3y^2+24z^2$ and hence $24n+4=(2x)^2+3(2y)^2+96z^2$.
 Since $w=(2x)^2+3(2y)^2\eq 4\ (\mo\ 8)$, by Lemma 3.2 we can write $w$ as $u^2+3v^2$ with $u$ and $v$ both odd. Now
 $24n+4=u^2+3v^2+96z^2$. Write $u$ or $-u$ in the form $6x-1$ with $x\in\Z$, and write $v$ or $-v$ in the form
 $4y-1$ with $y\in\Z$.
 Then we have
 $$24n+4=(6x-1)^2+3(4y-1)^2+96z^2$$
 and hence
 $$n=p_5(x)+p_6(y)+4p_4(z)=p_3(2y-1)+p_4(2z)+p_5(x).$$

(ii) By Dickson [D39, pp.\,112--113] or [JP] or [JKS], the quadratic form $6x^2+6y^2+z^2$ is regular
and it represents any nonnegative integer not in the set
$$\{8l+3:\ l\in\N\}\cup\{9^k(3l+2):\ k,l\in\N\}.$$
Thus there are $r,s,t\in\Z$ such that
$24n+1=6r^2+6s^2+t^2$. Clearly $2\nmid t$. As $6+6+1\not\eq1\ (\mo\ 8)$, we have $r\eq s\eq0\pmod2$.
Note that $\gcd(t,6)=1$. So, for some $x,y,z\in\Z$ we have
$$24n+1=24x^2+24y^2+(6z-1)^2,\ \ \t{i.e.,}\ n=p_4(x)+p_4(y)+p_5(z).$$

By Lemma 3.3, there are $u,v,w\in\Z$ with $w$ odd such that
$24n+1=6u^2+3v^2+w^2$. Clearly $2\mid v$. Since $3(2u^2+v^2)\eq1-w^2\eq0\ (\mo\ 8)$,
both $u$ and $v/2$ are even. Write $u=2x$ and $v/2=2y$ with $x,y\in\Z$.
As $\gcd(w,6)=1$, we can write $w$ or $-w$ in the form $6z-1$. Thus
$$24n+1=6u^2+3v^2+w^2=6(2x)^2+3(4y)^2+(6z-1)^2$$
and hence $$n=x^2+2y^2+p_5(z)=p_4(x)+2p_4(y)+p_5(z).$$

(iii)  (a) $m\in\N$ is said to be {\it $1$-good} (resp., {\it $3$-good})
if $m$ can be written in the form $x^2+3y^2+3z^2$ (resp., the form $3x^2+y^2+z^2$) with $x$ odd.
 Below we use induction to show that for $n=0,1,2,\ldots$ the number $12n+4$ is $1$-good and both $12n+4$ and $12n+8$
 are $3$-good.

 Clearly $4=1^2+3\times1^2$ is both $1$-good and $3$-good, and $8=3\times1^2+2^2+1^2$ is $3$-good.

 Now let $n\in\Z^+$ and assume that for any $m\in\{0,\ldots,n-1\}$ the number $12m+4$ is $1$-good
 and both $12m+4$ and $12m+8$ are $3$-good. We want to show that $12n+4$ is $1$-good and $12n+4r$
 is $3$-good, where $r\in\{1,2\}$.

 By Dickson [D39, pp.\,112--113], any nonnegative integer not of the form $3^{2k}(3l+2)$ (resp. $3^{2k+1}(3l+2)$)
 ($k,l\in\N$) can be written
 in the form $x^2+3y^2+3z^2$ (resp. $3x^2+y^2+z^2$) with $x,y,z\in\Z$.
 So, there are integers $x_0,y_0,z_0$ such that $12n+4=x_0^2+3y_0^2+3z_0^2$, also
  $12n+4r=3x_1^2+y_1^2+z_1^2$ for some $x_1,y_1,z_1\in\Z$.

 Let $i\in\{0,1\}$. If $x_i$ is even, then $y_i$ and $z_i$ are also even since $4\mid (y_i^2+z_i^2)$.
 If $x_0^2+3y_0^2\eq 4\ (\mo\ 8)$ (resp., $x_0^2+3z_0^2\eq 4\ (\mo\ 8)$), then by Lemma 3.2 we can write $x_0^2+3y_0^2$
 (resp. $x_0^2+3z_0^2$) in the form $a_0^2+3b_0^2$
 with $a_0$ and $b_0$ odd. If $3x_1^2+y_1^2\eq 4\ (\mo\ 8)$ (resp., $3x_1^2+z_1^2\eq 4\ (\mo\ 8)$),
 then by Lemma 3.2 we can write $3x_1^2+y_1^2$
 (resp. $3x_1^2+z_1^2$) in the form $3a_1^2+b_1^2$
 with $a_1$ and $b_1$ odd.

 Now it suffices to consider the case $x_i\eq y_i\eq z_i\eq0\pmod2$ and $x_i/2\eq y_i/2\eq z_i/2\pmod2$.
 Note that
 $$3n+1=\l(\f {x_0}2\r)^2+3\l(\f {y_0}2\r)^2+3\l(\f {z_0}2\r)^2,\ 3n+r=3\l(\f {x_1}2\r)^2+\l(\f {y_1}2\r)^2+\l(\f {z_1}2\r)^2.$$

 {\it Case} 1. $x_i/2,y_i/2,z_i/2$ are all odd.

 Without loss of generality, we may assume that $x_i/2\eq y_i/2\ (\mo\ 4)$
(otherwise we may use $-x_i/2$ instead of $x_i/2$).

Suppose $i=0$. Then $n$ is even since $3n+1$ is odd. Clearly,
$$12n+4=\l(\f{x_0+3y_0}2\r)^2+3\l(\f{x_0-y_0}2\r)^2+3z_0^2.$$
As $((x_0+3y_0)/2)^2+3z_0^2\eq 4-3((x_0-y_0)/2)^2\eq4\ (\mo\ 8)$,
by Lemma 3.2 we can write $((x_0+3y_0)/2)^2+3z_0^2=u_0^2+3v_0^2$ with $u_0$ and $v_0$ both odd. So
$12n+4=u_0^2+3v_0^2+3((x_0-y_0)/2)^2$ is $1$-good.

Now assume $i=1$. Then $3n+r$ is odd. Also,
$$12n+4r=4(3n+r)=3\l(\f{x_1-y_1}2\r)^2+\l(2x_1+\f{y_1-x_1}2\r)^2+z_1^2.$$
As $3((x_1-y_1)/2)^2+z_1^2\eq 4-(2x_1+(y_1-x_1)/2)^2\eq4\ (\mo\ 8)$,
by Lemma 3.2 we can write $3((x_1-y_1)/2)^2+z_1^2=3u_1^2+v_1^2$ with $u_1$ and $v_1$ both odd. So
$12n+4r=3u_1^2+v_1^2+((3x_1+y_1)/2)^2$ is $3$-good.

\medskip

{\it Case} 2. $x_i/2,y_i/2,z_i/2$ are all even.

Assume $i=1$. Then $3n+r\eq0\ (\mo\ 4)$ and hence $n=4m+r$ for some $m\in\N$.
 Since $3n+r=12m+4r<12n+4r$,
by the induction hypothesis there are $u,v,w\in\Z$ with $2\nmid u$ such that $3n+r=12m+4r=3u^2+v^2+w^2$.
Clearly $v\not\eq w\ (\mo\ 2)$. Without loss of generality, we assume that $w$ is odd.
Since $12n+4r=48m+16r=3(2u)^2+(2v)^2+(2w)^2$, we have $3(2u)^2+(2v)^2\eq -4w^2\eq4\ (\mo\ 8)$.
By Lemma 3.2, we can rewrite $3(2u)^2+(2v)^2$ as $3a^2+b^2$ with $a,b$ odd.
Thus $12n+4r=3a^2+b^2+(2w)^2$ is $3$-good.

By a similar argument, in the case $i=0$ we get that $12n+4$ is $1$-good.
\medskip

(b) Let $n\in\N$. By the Gauss-Legendre theorem,  $12n+5=r^2+s^2+t^2$ for some $r,s,t\in\Z$ with $s\eq t\ (\mo\ 2)$.
Note that $12n+5=r^2+2u^2+2v^2$ where $u=(s+t)/2$ and $v=(s-t)/2$. If $u\eq v\eq0\ (\mo\ 3)$
then $r^2\eq 5\eq 2\ (\mo\ 3)$
which is impossible. Without loss of generality we assume that $3\nmid v$.  Then $r^2+2u^2\eq 5-2v^2\eq0\ (\mo\ 3)$.
Since $r^2+2u^2>0$, by Lemma 2.1 we can write $r^2+2u^2$ in the form $a^2+2b^2$ with $a$ and $b$ relatively prime to 3.
Now we have $12n+5=a^2+2b^2+2v^2$ with $a,b,v$ relatively prime to 3. As $2(b^2+v^2)\eq5-a^2\eq 0\ (\mo\ 4)$,
we have $b\eq v\ (\mo\ 2)$. Choose $c\in\{v,-v\}$ such that $c\eq b\ (\mo\ 3)$. Then
$b-c=6z$ for some $z\in\Z$. Set $y=b+c$. Then
$$12n+5=a^2+2b^2+2c^2=a^2+(b+c)^2+(b-c)^2=a^2+y^2+(6z)^2.$$

(c) Fix $n\in\N$. By (b) there are $r,s,t\in\Z$ such that $12n+5=r^2+s^2+36t^2$. Clearly $r\not\eq s\ (\mo\ 2)$.
Note that $24n+10=(r+s)^2+(r-s)^2+72t^2$. Since $(r+s)^2+(r-s)^2\eq10\eq1\ (\mo\ 3)$, one of $r+s$ and $r-s$
is congruent to 3 mod 6, and the other is relatively prime to 6. So there are $u,v\in\Z$ such that
$$24n+10=(3(2u+1))^2+(6v-1)^2+72t^2,\ \t{i.e.,}\ n=3p_3(u)+3p_4(t)+p_5(v).$$

Choose $x,y,z\in\Z$ with $36n+29=x^2+y^2+(6z)^2$. Then $x\not\eq y\ (\mo\ 2)$ and
$72n+58=(x+y)^2+(x-y)^2+72z^2$.
Since $58\eq1\ (\mo\ 3)$, without loss of generality we assume that $x+y\eq0\ (\mo\ 3)$ and $x-y\eq1\ (\mo\ 3)$.
Note that $(x-y)^2\eq 58\eq 7^2\ (\mo\ 9)$ and hence $x-y\eq \pm 7\ (\mo\ 18)$. So, there are $x_0,y_0\in\Z$ such that
$$72n+58=(3(2x_0+1))^2+(18y_0-7)^2+72z^2$$
and hence $n=p_3(x_0)+p_4(z)+p_{11}(y_0)$.
\medskip

(iv) Assume that $20n+r$ is not a square.
By the Gauss-Legendre theorem, there are $x,y,z\in\Z$ such that $20n+r=x^2+y^2+z^2$. It is easy to see that
we cannot have $x^2,y^2,z^2\not\eq r\ (\mo\ 5)$. Without loss of generality we assume that $x^2\eq r\ (\mo\ 5)$.
If $2\nmid yz$, then $x^2+y^2+z^2\eq x^2+2\not\eq r\ (\mo\ 4)$.  So $y$ or $z$ is even.
Without loss of generality we assume that $2\mid z$. Since $20n+r\not=x^2$, by Lemma 3.5 we can write
$y^2+z^2=y^2+4(z/2)^2$ in the form  $s^2+4t^2$ with $s,t\in\Z$ not all divisible by 5.
As $s^2+4t^2=y^2+z^2\eq0\ (\mo\ 5)$,
we have $20n+r=x^2+s^2+(2t)^2$ with $5\nmid xst$.
Therefore there are integers $u,v,w$ relatively prime to 5
such that $20n+r=(2u)^2+(2v)^2+w^2$. If $(2u)^2,(2v)^2\not\eq r\ (\mo\ 5)$, then
$w^2\eq r-(-r-r)\ (\mo\ 5)$ which is impossible. Without loss of generality, we assume that $(2u)^2\eq r\ (\mo\ 5)$.
Then $v^2\eq -(2v)^2\eq w^2\ (\mo\ 5)$. Simply let $v\eq w\ (\mo\ 5)$ (otherwise we may change the sign of $w$).
Set $a=(4v+w)/5$ and $b=(w-v)/5$. Then $a-b=v$ and $a+4b=w$. Thus
$$20n+r=4u^2+w^2+4v^2=4u^2+(a+4b)^2+4(a-b)^2=4u^2+5a^2+5(2b)^2.$$

Now suppose that $20(n+1)+1$ is not a square. By the above, there are $u,v,w\in\Z$ such that $20n+21=5u^2+5v^2+4w^2$.
Clearly $u\not\eq v\ (\mo\ 2)$, $4w^2\eq21\eq 4^2\ (\mo\ 5)$ and hence $w\eq\pm 2\ (\mo\ 5)$. Write $u+v=2x+1$,
$u-v=2y+1$, and $w$ or $-w$ in the form $5z-2$, where $x,y,z$ are integers. Then
$$40n+42=5(u+v)^2+5(u-v)^2+8w^2=5(2x+1)^2+5(2y+1)^2+8(5z-2)^2$$
and hence $n=p_3(x)+p_3(y)+p_{12}(z)$.
\medskip

(v) Let $n\in\N$. Observe that
$$n=p_3(x)+p_4(y)+p_9(z)\iff 56n+32=7(2x+1)^2+56y^2+(14z-5)^2.$$
If $56n+32=7x^2+56y^2+z^2$ for some $x,y,z\in\Z$, then $7x^2+z^2>0$ and $8\mid (7x^2+z^2)$, hence
by Lemma 3.6 we can write $7x^2+z^2$ as $7a^2+b^2$ with $a$ and $b$ both odd, and thus
$b\eq\pm 5\ (\mo\ 14)$ since
$b^2\eq z^2\eq 32\eq 5^2\ (\mo\ 7)$. Therefore
$$\align &n=p_3(x)+p_4(y)+p_9(z)\ \t{for some}\ x,y,z\in\Z
\\\iff&56n+32=x^2+7y^2+56z^2\ \t{for some}\ x,y,z\in\Z.
\endalign$$

 Now assume that $7n+4$ is not squarefree. By [WP, Theorem 3], there are $u,v,w\in\Z$ such that
$7n+4=u^2+7v^2+7w^2$ and hence
$$56n+32=8(7n+4)=8(u^2+7v^2)+56w^2=(u+7v)^2+7(u-v)^2+56w^2$$
as desired.

\medskip

(vi) Let $n\in\N$. Clearly
$$n=p_3(x)+p_4(y)+p_{18}(z)\iff32n+53=(2(2x+1))^2+32y^2+(16z-7)^2.$$
On the other hand, if $32n+53=x^2+y^2+32z^2$ with $x,y,z\in\Z$, then we can write $\{x,y\}=\{2u,2v+1\}$ with $u,v\in\Z$,
and we have $2\nmid u$ (since $4u^2\eq 53-(2v+1)^2\eq4\ (\mo\ 8)$) and $2v+1\eq\pm 7\ (\mo\ 16)$
(since $(2v+1)^2\eq 53-4u^2\eq 53-4=7^2\ (\mo\ 32)$). So
$$\align&n=p_3(x)+p_4(y)+p_{18}(z)\ \t{for some}\ x,y,z\in\Z
\\\iff&32n+53=\f{(x+y)^2+(x-y)^2}2+32z^2\ \t{for some}\ x,y,z\in\Z
\\\iff&8(8n+13)+2=64n+106=x^2+y^2+64z^2\ \t{for some}\ x,y,z\in\Z.
\endalign$$

Observe that
$$\align &n=p_3(x)+p_4(y)+p_{27}(z)
\\\iff&200n+554=(5(2x+1))^2+200y^2+(50z-23)^2.
\endalign$$
On the other hand, if $200n+554=x^2+y^2+200z^2$, then $x^2+y^2\eq554\eq2\ (\mo\ 4)$ and hence $x\eq y\eq1\ (\mo\ 2)$,
also $x^2+y^2\eq 4\ (\mo\ 25)$ and hence we can write $\{x,y\}=\{5u,v\}$ with $u,v\in\Z$ and $v^2\eq 4\eq(-23)^2\ (\mo\ 25)$. Thus
$$\align&n=p_3(x)+p_4(y)+p_{27}(z)\ \t{for some}\ x,y,z\in\Z
\\\iff&200n+554=(x+y)^2+(x-y)^2+200z^2\ \t{for some}\ x,y,z\in\Z
\\\iff&4(25n+69)+1=x^2+y^2+100z^2\ \t{for some}\ x,y,z\in\Z.
\endalign$$

 By [KS, Corollary 1.9(ii)],
the sets
$$\align S=&\{n\in\N:\ n\not=p_3(x)+p_3(y)+8z^2\ \t{for any}\ x,y,z\in\Z\}
\\=&\{n\in\N:\ 8n+2\not= x^2+y^2+64z^2\ \t{for any}\ x,y,z\in\Z\}
\endalign$$
and
$$\align T=&\{n\in\N:\ n\not=x^2+2p_3(y)+25z^2\ \t{for any}\ x,y,z\in\Z\}
\\=&\{n\in\N:\ 4n+1\not=x^2+y^2+100z^2\ \t{for any}\ x,y,z\in\Z\}
\endalign$$
are finite. Our computation suggests that
$S=\{5,\,40,\,217\}$ and
$$\align T=&\{5,\,8,\,14,\,17,\,19,\,23,\,33,\,44,\,75,\,77,\,96,\,147,\,180,\,195,
\\&\ 203,\,204,\,209,\,222,\,482,\,485,\,495,\,558,\,720,\,854,\,1175\}.
\endalign$$
As mentioned in [K09, Section 3] or at the end of [KS, Example 1.1] , under the GRH one can use the argument of Ono
and Soundararajan [OnoS] to show that $S$ and $T$ indeed only contain the listed elements.
So part (v) of Theorem 1.7 follows.

\medskip

 In view of the above, we have completed the proof of Theorem 1.7.
 \qed

 \medskip

\heading{4. Proof of Theorem 1.14}\endheading

\proclaim{Lemma 4.1} Let $n\eq 2\ (\mo\ 3)$ be a nonnegative integer not of the form $4^k(8l+7)$ with $k,l\in\N$.
 Then there are $x,y,z\in\Z$ such that
 $$2n=x^2+9(y^2+2z^2)=x^2+3(y+2z)^2+6(y-z)^2.$$
 \endproclaim
 \Proof. By the Gauss-Legendre theorem, there are $u,v,w\in\Z$ such that $n=u^2+v^2+w^2$.
 Since $n\eq2\ (\mo\ 3)$, exactly one of $u,v,w$ is divisible by 3. Without loss of generality, we assume that
 $w=3z$ with $z\in\Z$. Write $\{u+v,u-v\}=\{x,3y\}$ with $x,y\in\Z$. Then
 $2n=(u+v)^2+(u-v)^2+2w^2=x^2+9y^2+18z^2$ and hence the desired result follows. \qed
\medskip

For $a,b,c\in\Z^+$ we define
$$E(ax^2+by^2+cz^2)=\{n\in\N:\ n\not=ax^2+by^2+cz^2\ \t{for any}\ x,y,z\in\Z\}.$$
By Dickson [D27] and [D39, pp.\,112--113], we have
 $$\align E(x^2+3y^2+18z^2)=&\bigcup_{k,l\in\N}\{3l+2,9l+6,4^k(16l+10)\},\tag4.1
\\E(x^2+3y^2+2z^2)=&\{4^k(16l+10):\ k,l\in\N\},\tag4.2
\\E(6x^2+6y^2+z^2)=&\{8l+3:\ l\in\N\}\cup\{9^k(3l+2):\ k,l\in\N\},\tag4.3
\\E(3x^2+3y^2+z^2)=&\{9^k(3l+2):\ k,l\in\N\},\tag4.4
\\E(3x^2+12y^2+2z^2)=&\{16l+6:\ l\in\N\}\cup \{9^k(3l+1):\ k,l\in\N\},\tag4.5
\\E(2x^2+3y^2+6z^2)=&\{3l+1:\ l\in\N\}\cup \{4^k(8l+7):\ k,l\in\N\},\tag4.6
\\E(2x^2+3y^2+48z^2)=&\bigcup_{k,l\in\N}\{8l+1,8l+7,16l+6,16l+10,64l+24,9^k(3l+1)\},\tag4.7
\\E(6x^2+18y^2+z^2)=&\bigcup_{l\in\N}\{3l+2,\, 9l+3\}\cup\{4^k(8l+5): k,l\in\N\},\tag4.8
\\E(9x^2+24y^2+z^2)=&\bigcup_{l\in\N}\{3l+2,4l+3,8l+6\}\cup\{9^k(9l+3): k,l\in\N\},\tag4.9
\\E(2x^2+2y^2+3z^2)=&\{8l+1:\ l\in\N\}\cup\{9^k(9l+6): k,l\in\N\},\tag4.10
\\ E(5x^2+5y^2+z^2)=&\bigcup_{l\in\N}\{5l+2,\, 5l+3\}\cup\{4^k(8l+7): k,l\in\N\},\tag4.11
\\E(5x^2+10y^2+2z^2)=&\{8l+3:\,l\in\N\}\cup\bigcup_{k,l\in\N}\{25^k(5l+1),25^k(5l+4)\},\tag4.12
\\E(5x^2+40y^2+z^2)=&\bigcup_{k,l\in\N}\{4l+3,\, 8l+2,25^k(5l+2),25^k(5l+3)\},\tag4.13
\\E(5x^2+15y^2+6z^2)=&\bigcup_{k,l\in\N}\{9^k(3l+1),\,4^k(16l+14),\,25^k(5l+2),\,25^k(5l+3)\},\tag4.14
\\E(3x^2+3y^2+4z^2)=&\{4l+1:\ l\in\N\}\cup\{9^k(3l+2):\, k,l\in\N\},\tag4.15
\endalign$$
and
$$\align
E(3x^2+24y^2+8z^2)=&\bigcup_{k,l\in\N}\{3l+1,4l+1,4l+2,4^k(8l+7)\},\tag4.16
\\E(3x^2+48y^2+8z^2)=&\bigcup_{k,l\in\N}\{4l+1,4l+2,8l+7,64l+24,9^k(3l+1)\},\tag4.17
\\E(3x^2+72y^2+8z^2)=&\bigcup_{k,l\in\N}\{3l+1,4l+1,4l+2,8l+7,32l+4,9^k(9l+6)\},\tag4.18
\\E(3x^2+36y^2+4z^2)=&\bigcup_{k,l\in\N}\{3l+2,4l+1,4l+2,9^k(9l+6)\},\tag4.19
\\E(2x^2+4y^2+z^2)=&\{4^k(16l+14):\ k,l\in\N\},\tag4.20
\\E(10x^2+2y^2+5z^2)=&\bigcup_{k,l\in\N}\{8l+3,\, 25^k(5l+1),\, 25^k(5l+4)\},\tag4.21
\\E(5x^2+40y^2+8z^2)=&\bigcup_{k,l\in\N}\{4l+2,4l+3,8l+1,32l+12,5^{2k}(5l+1),5^{2k}(5l+4)\}.\tag4.22
\endalign$$

\medskip
\noindent{\it Proof of Theorem 1.14}. (i) Clearly,
$$\align &n=p_3(x)+p_3(y)+2p_5(z)\ \t{for some}\ x,y,z\in\Z
\\\iff&24n+8=3(2x+1)^2+3(2y+1)^2+2(6z-1)^2\ \t{for some}\ x,y,z\in\Z
\\\iff&24n+8=3x^2+3y^2+2z^2\ \ \t{for some}\ x,y,z\in\Z\ \t{with}\ 2\nmid z
\\\iff&24n+8=3(x+y)^2+3(x-y)^2+2z^2\ \ \t{for some}\ x,y,z\in\Z\ \t{with}\ 2\nmid z
\\\iff&12n+4=3x^2+3y^2+z^2\ \ \t{for some}\ x,y,z\in\Z\ \t{with}\ 2\nmid z.
\endalign$$
Also,
$$\align &n=p_3(x)+2p_4(y)+p_5(z)\ \t{for some}\ x,y,z\in\Z
\\\iff&24n+4=3(2x+1)^2+48y^2+(6z-1)^2\ \t{for some}\ x,y,z\in\Z
\\\iff&24n+4=3x^2+12y^2+z^2\ \ \t{for some}\ x,y,z\in\Z\ \t{with}\ 2\nmid z.
\endalign$$
(If $24n+4=3x^2+12y^2+z^2$ with $2\nmid z$, then $\gcd(z,6)=1$, and also $2\mid y$
since $3x^2+z^2\eq3+1=4\ (\mo\ 8)$.) When $24n+4=3x^2+3y^2+z^2$ with $z$ odd, one of $x$ and $y$ is even.
Therefore, by the above and Theorem 1.7(iii), both $p_3+p_3+2p_5$ and $p_3+2p_4+p_5$ are universal over $\Z$.

 By Theorem 1.7(iii), there are $u,v,w\in\Z$ with $2\nmid w$ such that $12n+8=u^2+v^2+3w^2$.
Clearly $u\not\eq v\ (\mo\ 2)$ and $24n+16=(u+v)^2+(u-v)^2+6w^2$. Since $(u+v)^2+(u-v)^2\eq 16\eq1\ (\mo\ 3)$,
exactly one of $u+v$ and $u-v$ is divisible by 3. Thus there are $x,y,z\in\Z$ such that
$$24n+16=6(2x+1)^2+(3(2y+1))^2+(6z-1)^2$$
and hence $n=2p_3(x)+3p_3(y)+p_5(z)$.

 Let $m\in\{1,2,3\}$. As $24n+6m+1\not\eq 1\ (\mo\ 8)$ is not a square, by Lemma 3.3 there are $u,v,w\in\Z$ such that
 $24n+6m+1=(2u)^2+3v^2+6w^2$. Clearly $2\nmid v$.
 Also, $w\not\eq m\ (\mo\ 2)$ since $6w^2\eq 6m+1-3v^2\eq 2(m-1)\ (\mo\ 4)$.
 Note that
 $$4u^2\eq 6m+1-3v^2-6w^2\eq -2m-2+2(m-1)^2\eq 2m(m+1)\ (\mo\ 8).$$
 Hence $2\nmid u$ if $m\in\{1,2\}$, and $2\mid u$ if $m=3$. Clearly $3\nmid u$.
 In the case $m=1$, there are $x,y,z\in\Z$ such that
 $$24n+7=4u^2+3v^2+6w^2=6(2x)^2+3(2y+1)^2+4(6z-1)^2$$
 and hence $n=x^2+p_3(y)+4p_5(z)$. Though $p_3+p_4+4p_5$ is universal over $\Z$,
 it is not universal over $\N$. When $m=2$, there are $x,y,z\in\Z$ such that
 $$24n+13=4u^2+3v^2+6w^2=3(2x+1)^2+6(2y+1)^2+4(6z-1)^2$$
 and hence $n=p_3(x)+2p_3(y)+4p_5(z)$.
 In the case $m=3$, both $u$ and $w$ are even, hence there are $x,y,z\in\Z$ such that
 $$24n+19=4u^2+3v^2+6w^2=3(2x+1)^2+6(2y)^2+4(2(3z-1))^2$$
 and thus $n=p_3(x)+p_4(y)+2p_8(z)$.

 By Theorem 1.7(iii), there are $u,v,w\in\Z$ with $2\nmid u$ such that $24n+16=u^2+3v^2+3w^2$.
 Clearly $v\not\eq w\ (\mo\ 2)$. Without loss of generality we assume that $2\nmid v$ and $2\mid w$.
 As $3w^2\eq 16-u^2-3v^2\eq 4\ (\mo\ 8)$, $w/2$ is odd. Note that $u$ is relatively prime to 6.
 Thus there are $x,y,z\in\Z$ such that
$$24n+16=3(2x+1)^2+12(2y+1)^2+(6z-1)^2$$
and hence $n=p_3(x)+4p_3(y)+p_5(z)$.

By (4.1), there are $u,v,w\in\Z$ such that $24n+22=u^2+3v^2+18w^2$. Since $u\eq v\ (\mo\ 2)$,
 we have $u^2+3v^2\eq0\ (\mo\ 4)$ and hence $w$ is odd. As $u^2+3v^2\eq 22-18w^2\eq 4\ (\mo\ 8)$,
by Lemma 3.2 there are odd numbers $a$ and $b$ such that $u^2+3v^2=a^2+3b^2$.
Clearly $\gcd(a,6)=1$. So there are $x,y,z\in\Z$ such that
$$24n+22=3b^2+18w^2+a^2=3(2x+1)^2+18(2y+1)^2+(6z-1)^2$$
and hence $n=p_3(x)+6p_3(y)+p_5(z)$.

By (4.2), for $m\in\{0,1,2\}$, there are $u,v,w\in\Z$ such that $24n+8m+4=u^2+3v^2+2w^2$.
Clearly $4\mid (u^2+3v^2)$ and hence $2\mid w$. So $u^2+3v^2\eq 4\ (\mo\ 8)$. By Lemma 3.2,
$u^2+3v^2=a^2+3b^2$ for some odd integers $a$ and $b$. By Lemma 2.1,
$a^2+2w^2$ can be written as $s^2+2t^2$ with $s$ or $t$ not divisible by 3.
Since $24n+8m+4=a^2+3b^2+2w^2=s^2+3b^2+2t^2$, we have $4\mid (s^2+3b^2)$ and hence $2\mid t$.
In the case $m=0$, we have $s^2+2t^2\eq 4\eq1\ (\mo\ 3)$ and hence $s\not\eq t\eq0\ (\mo\ 3)$, thus
there are $x,z\in\Z$ such that
 $$24n+4=s^2+3b^2+2\times 6^2\l(\f t6\r)^2=3(2x+1)^2+72y^2+(6z-1)^2$$
 and hence $n=p_3(x)+3p_4(y)+p_5(z)$.
If $m=1$, then both $s$ and $t$ are relatively prime to 3, so there are $x,y,z\in\Z$ such that
$$24n+12=s^2+3b^2+2t^2=3(2x+1)^2+(6y-1)^2+8(3z-1)^2$$
and hence $n=p_3(x)+p_5(y)+p_8(z)$.
When $m=2$, we have $s^2+2t^2\eq 20\eq2\ (\mo\ 3)$ and hence $t\not\eq s\eq0\ (\mo\ 3)$, thus
there are $x,z\in\Z$ such that
 $$24n+20=s^2+3b^2+2t^2=(3(2x+1))^2+3(2y+1)^2+2(2(3z-1))^2$$
 and hence $n=3p_3(x)+p_3(y)+p_8(z)$. Though $p_3+3p_3+p_8$ is universal over $\Z$, it is not universal over $\N$.

(ii) By (4.3),
for some $u,v,w\in\Z$ we have
$$24n+7=6u^2+6v^2+w^2=3(u+v)^2+3(u-v)^2+w^2.$$
It is easy to see that $\gcd(w,6)=1$ and $u+v\eq u-v\eq1\pmod2$. Thus, there are $x,y,z\in\Z$ such that
$$24n+7=3(2x+1)^2+3(2y+1)^2+(6z-1)^2$$
and hence $n=p_3(x)+p_3(y)+p_5(z)$.

By Lemma 4.1, there are $u,v,w\in\Z$ such that $2(12n+5)=u^2+3v^2+6w^2$.
It is easy to see that $2\nmid w$ and $u^2+3v^2\eq4\ (\mo\ 8)$. By Lemma 3.2, $u^2+3v^2=a^2+3b^2$ for some odd integers $a$ and $b$.
So there are $x,y,z\in\Z$ such that
$$24n+10=a^2+3b^2+6w^2=(6z-1)^2+3(2x+1)^2+6(2y+1)^2$$
and hence $n=p_3(x)+2p_3(y)+p_5(z)$.

In light of (4.4), $24n+10=3u^2+3v^2+w^2$ for some $u,v,w\in\Z$. It is easy to see that $2\mid w$ and $2\nmid uv$.
Since $w^2\eq10-3(u^2+v^2)\eq 4\ (\mo\ 8)$, $w/2$ is odd and hence $\gcd(w/2,6)=1$.
Thus, for some $x,y,z\in\Z$ we have
$$24n+10=3(2x+1)^2+3(2y+1)^2+4(6z-1)^2$$
and hence  $n=p_3(x)+p_3(y)+4p_5(z)$.

By (4.5), $24n+17=3u^2+12v^2+2w^2$ for some $u,v,w\in\Z$.
It is easy to see that $2\nmid uv$ and $\gcd(w,6)=1$. So, there are $x,y,z\in\Z$ such that
$$24n+17=3(2x+1)^2+12(2y+1)^2+2(6z-1)^2$$
and hence
$n=p_3(x)+4p_3(y)+2p_5(z)$.

In view of (4.6), $24n+5=2u^2+3v^2+6w^2$ for some $u,v,w\in\Z$. The equality modulo $8$ yields that $2\nmid uv$ and $2\mid w$.
Note that $\gcd(u,6)=1$. So, for some $x,y,z\in\Z$
we have
$$24n+5=2(6z-1)^2+3(2x+1)^2+24y^2$$
and hence $n=p_3(x)+p_4(y)+2p_5(z)$.

By (4.7), $24n+5=2u^2+3v^2+48w^2$ for some $u,v,w\in\Z$. Clearly $2\nmid v$,
and $\gcd(u,6)=1$.  So, for some $x,z\in\Z$
we have
$$24n+5=2(6z-1)^2+3(2x+1)^2+48w^2$$
and hence $n=p_3(x)+2p_4(w)+2p_5(z)$.

In light of (4.8), there are $u,v,w\in\Z$ such that
$24n+7=6u^2+18v^2+w^2$. Clearly $\gcd(w,6)=1$, $2\nmid u$ and $2\mid v$. Thus, for some $x,y,z\in\Z$ we have
$$24n+7=6(2x+1)^2+18(2y)^2+(6z-1)^2$$
and hence $n=2p_3(x)+3p_4(y)+p_5(z).$

By (4.9), there are $u,v,w\in\Z$ such that
$24n+10=9u^2+24v^2+w^2$. Clearly $u^2+w^2\eq 9u^2+w^2\eq 10\eq 2\ (\mo\ 8)$.
So $2\nmid u$ and $\gcd(w,6)=1$. Therefore, for some $x,z\in\Z$ we have
$$24n+10=9(2x+1)^2+24v^2+(6z-1)^2$$
and hence $n=3p_3(x)+p_4(v)+p_5(z)$.

In view of (4.10), for any given $\da\in\{0,1\}$, there are $u,v,w\in\Z$ such that
$$24n+8\da+5=2u^2+2v^2+3w^2=(u+v)^2+(u-v)^2+3w^2.$$  Clearly $w$ and $u\pm v$ are all odd. Note that $x^2\eq5\ (\mo\ 3)$
for no $x\in\Z$. In the case $\da=0$, we have  $\gcd(u\pm v,6)=1$, so there are $x,y,z\in\Z$ such that
$$24n+5=(6y-1)^2+(6z-1)^2+3(2x+1)^2$$
and hence $n=p_3(x)+p_5(y)+p_5(z).$
In the case $\da=1$, we have $(u+v)^2+(u-v)^2\eq13\eq1\ (\mo\ 3)$, hence exactly one of $u+v$ and $u-v$
is divisible by 3, thus there are $x,y,z\in\Z$ such that
$$24n+13=3(2x+1)^2+9(2y+1)^2+(6z-1)^2$$
and hence $n=p_3(x)+3p_3(y)+p_5(z).$

(iii) Observe that $p_7(z)=(5z^2-3z)/2$ and $40p_7(z)+9=(10z-3)^2$.

By (4.11), there are $u,v,w\in\Z$ such that
$40n+19=5u^2+5v^2+w^2$. Since $u^2+v^2+w^2\eq3\pmod 4$, we have $2\nmid uvw$.
Clearly $w^2\eq 19\eq 3^2\ (\mo\ 5)$, thus $w$ or $-w$ is congruent to
$-3$ mod $10$.
So, for some $x,y,z\in\Z$ we have
$$40n+19=5(2x+1)^2+5(2y+1)^2+(10z-3)^2$$
and hence $n=p_3(x)+p_3(y)+p_7(z).$

In light of (4.12), there are $u,v,w\in\Z$ such that
$40n+33=5u^2+10v^2+2w^2$. This equality modulo $8$ yields that $2\nmid uvw$. Note that $2w^2\eq33\eq18\ (\mo\ 5)$ and hence $w\eq\pm3\ (\mo\ 5)$.
So there are $x,y,z\in\Z$ such that
$$40n+33=5(2x+1)^2+10(2y+1)^2+2(10z-3)^2$$
and hence $n=p_3(x)+2p_3(y)+2p_7(z)$.

By (4.13), there are $u,v,w\in\Z$ such that
$40n+14=5u^2+40v^2+w^2$. Clearly $u\eq w\eq1\ (\mo\ 2)$ since $4\nmid14$.
Note that $w^2\eq14\eq3^2\ (\mo\ 5)$ and hence $w\eq\pm3\ (\mo\ 10)$.
So, for some $x,z\in\Z$ we have
$$40n+14=5(2x+1)^2+40v^2+(10z-3)^2$$
and hence $n=p_3(x)+p_4(v)+p_7(z)$.

In view of (4.14), there are $u,v,w\in\Z$ such that $120n+74=5u^2+15v^2+6w^2$. As $4\mid (u^2+3v^2)$, we have $6w^2\eq74\ (\mo\ 4)$
 and hence $2\nmid w$. Note that $5(u^2+3v^2)\eq74-6w^2\eq 68\eq20\ (\mo\ 8)$ and hence $u^2+3v^2\eq4\ (\mo\ 8)$.
 By Lemma 3.2, $u^2+3v^2=a^2+3b^2$ for some odd integers $a$ and $b$. Now, $120n+74=5a^2+15b^2+6w^2$.
Clearly $\gcd(a,6)=1$. Also, $w^2\eq 6w^2\eq 74\eq 3^2\ (\mo\ 5)$ and hence $w\eq\pm 3\ (\mo\ 5)$.
So, for some $x,y,z\in\Z$ we have
$$120n+74=15(2x+1)^2+5(6y-1)^2+6(10z-3)^2$$
and hence $n=p_3(x)+p_5(y)+2p_7(z).$

(iv) As $p_8(z)=3z^2-2z$, we have $3p_8(z)+1=(3z-1)^2.$

In light of (4.15), there are $u,v,w\in\Z$ such that
$12n+7=3u^2+3v^2+4w^2$. It follows that
$$24n+14=6u^2+6v^2+8w^2=3(u+v)^2+3(u-v)^2+8w^2.$$
Clearly, $u\pm v\eq1\ (\mo\ 2)$ and $\gcd(w,3)=1$.
Thus, for some $x,y,z\in\Z$ we have
$$24n+14=3(2x+1)^2+3(2y+1)^2+8(3z-1)^2$$
and hence $n=p_3(x)+p_3(y)+p_8(z)$.

By (4.6) there are $u,v,w\in\Z$ such that $24n+17=3u^2+6v^2+2w^2$.
Clearly $2\nmid u$. Since $2(3v^2+w^2)\eq17-3u^2\eq 6\ (\mo\ 8)$, we
have $2\nmid v$ and $2\mid w$. Note that $\gcd(w/2,3)=1$. Thus, for
some $x,y,z\in\Z$ we have
$$24n+17=3(2x+1)^2+6(2y+1)^2+8(3z-1)^2$$
and hence $n=p_3(x)+2p_3(y)+p_8(z).$

Let $m\in\{1,2,3\}$. By (4.16)-(4.18), there are $u,v,w\in\Z$ such that
$24n+11=3u^2+24mv^2+8w^2$. Clearly $2\nmid u$ and $3\nmid w$. So, for some $x,y,z\in\Z$ we have
$$24n+11=3(2x+1)^2+24my^2+8(3z-1)^2$$
and hence $n=p_3(x)+mp_4(y)+p_8(z).$

In view of (4.19), there are $u,v,w\in\Z$ such that
$12n+7=3u^2+36v^2+4w^2$. Clearly $2\nmid u$ and $3\nmid w$.
So, for some $x,y,z\in\Z$ we have
$$12n+7=3(2x+1)^2+36y^2+4(3z-1)^2$$
and hence $n=2p_3(x)+3p_4(y)+p_8(z).$

(v) By (4.4) there are $a,b,c\in\Z$ such that $168n+103=a^2+3b^2+3c^2$. Since $a^2\not\eq103\eq3\ (\mo\ 4)$,
$b$ and $c$ cannot be all even. Without loss of generality we assume that $2\nmid c$.
Note that $a^2+3b^2\eq103-3c^2\eq 4\ (\mo\ 8)$. By Lemma 3.2, $a^2+3b^2=c^2+3d^2$ for some odd integers
$c$ and $d$. Without loss of generality we simply suppose that $a$ and $b$ are odd.

We claim that there are odd integers $a_0,b_0,c_0$ such that
$168n+103=a_0^2+3b_0^2+3c_0^2$ and $c_0\eq\pm 2\ (\mo\ 7)$.
Assume that $b,c\not\eq\pm2\ (\mo\ 7)$.
Since $a^2+3b^2+3c^2\eq 103\eq5\ (\mo\ 7)$, we cannot have $a,b,c\not\eq0\pmod 7$.
Without loss of generality, we suppose that $a$ or $b$ is divisible by 7,
and $a\eq b\ (\mo\ 4)$ (otherwise we use $-a$ instead of $a$).
Since $c\not\eq\pm2\ (\mo\ 7)$, we cannot have $a\eq b\eq0\ (\mo\ 7)$.
As $a^2+3b^2+3c^2\eq 5\ (\mo\ 7)$, the integer in $\{a,b\}$ not divisible by 7
must be congruent to 3 or $-3$ modulo 7.
Clearly, $(2a)^2+3(2b)^2=(a-3b)^2+3(a+b)^2$, $a-3b\eq a+b\eq 2\ (\mo\ 4)$ and
$a+b\eq \pm 4\ (\mo\ 7)$. Thus there are odd integers
$a_0,b_0$ such that $a^2+3b^2=a_0^2+3b_0^2$ and $b_0\eq\pm 2\ (\mo\ 7)$. This proves the claim.

By the above, there are odd integers $u,v,w\in\Z$ such that
$168n+103=u^2+3v^2+3w^2$  and $w\eq\pm 2\ (\mo\ 7)$. As $3w^2\eq 12\eq 103\ (\mo\ 7)$,
we have $u^2\eq-3v^2\eq (2v)^2\ (\mo\ 7)$. Without loss of generality we assume that $u\eq -2v\ (\mo\ 7)$.
Set $s=(u+2v)/7$ and $t=2s-v=(2u-3v)/7$. Then $2\nmid st$ and
$$u^2+3v^2=(3s+2t)^2+3(2s-t)^2=21s^2+7t^2.$$
Therefore $168n+103=21s^2+7t^2+3w^2$. Clearly $\gcd(t,6)=1$ and
$w\eq \pm 5\ (\mo\ 14)$. So there are $x,y,z\in\Z$ such that
$$168n+103=21(2x+1)^2+7(6y-1)^2+3(14z-5)^2$$
and hence $n=p_3(x)+p_5(y)+p_9(z)$.

(vi) Observe that $p_{10}(z)=4z^2-3z$ and $16p_{10}(z)+9=(8z-3)^2.$

By the Gauss-Legendre theorem, there are $u,v,w\in\Z$ such that
$16n+13=u^2+v^2+w^2$ with $w$ odd.
Note that $16n+13=2(s^2+t^2)+w^2$ where $s=(u+v)/2$ and $t=(u-v)/2$. As $2(s^2+t^2)\eq 13-w^2\eq12\eq4\ (\mo\ 8)$,
we have $s\eq t\eq1\ (\mo\ 2)$. Clearly $w^2\eq 13-2(s^2+t^2)\eq 3^2\ (\mo\ 16)$ and hence
$w\eq\pm 3\ (\mo\ 8)$. So, for some $x,y,z\in\Z$ we have
$$16n+13=2(2x+1)^2+2(2y+1)^2+(8z-3)^2$$
and hence $n=p_3(x)+p_3(y)+p_{10}(z).$

In view of (4.20), there are $u,v,w\in\Z$ such that
$16n+15=2u^2+4v^2+w^2$. This equality modulo 8 yields that $2\nmid uvw$. As $w^2\eq 15-2u^2-4v^2\eq 9\ (\mo\ 16)$, either $w$ or $-w$ is congruent to
$-3$ mod $8$.
So, for some $x,y,z\in\Z$ we have
$$16n+15=2(2x+1)^2+4(2y+1)^2+(8z-3)^2$$
and hence $n=p_3(x)+2p_3(y)+p_{10}(z)$.

By (4.20), $16n+11=2u^2+4v^2+w^2$ for some $u,v,w\in\Z$. This equality modulo 8 yields that $2\nmid uw$ and
$2\mid v$. Note that $w^2\eq11-2u^2\eq9\ (\mo\ 16)$ and hence $w\eq\pm 3\ (\mo\ 8)$.
Thus, for some $x,y,z\in\Z$ we have
$$16n+11=2(2x+1)^2+4(2y)^2+(8z-3)^2$$
and hence $n=p_3(x)+p_4(y)+p_{10}(z).$

In light of (4.21), there are $u,v,w\in\Z$ such that
$80n+73=10u^2+2v^2+5w^2$. This equality modulo 8 yields that $2\nmid uvw$. Note that
$$5w^2\eq73-10u^2-2v^2\eq73-12\eq5\times3^2\ (\mo\ 16)$$
and hence $w\eq\pm 3\ (\mo\ 8)$. Also, $2v^2\eq 73\eq 2\times3^2\ (\mo\ 5)$ and hence $v\eq\pm3\ (\mo\ 5)$.
So, for some $x,y,z\in\Z$ we have
$$80n+73=10(2x+1)^2+2(10y-3)^2+5(8z-3)^2$$
and hence $n=p_3(x)+p_7(y)+p_{10}(z).$

By (4.22), there are $u,v,w\in\Z$ such that
$40n+37=5u^2+40v^2+8w^2$.
Clearly $2\nmid u$. Also,
$8w^2\eq37\eq 32\ (\mo\ 5)$ and hence $w\eq\pm2\ (\mo\ 5)$.
So, for some $x,y,z\in\Z$ we have
$$40n+37=5(2x+1)^2+40y^2+8(5z-2)^2$$
and hence $n=p_3(x)+p_4(y)+p_{12}(z)$
since
$5p_{12}(z)=5(5z^2-4z)=(5z-2)^2-4.$

In view of the above, we have completed the proof of Theorem 1.14. \qed

\heading{5. Proofs of Theorems 1.17 and 1.20}\endheading

\medskip
\noindent{\it Proof of Theorem 1.17}. Let $n$ be any nonnegative integer.

(i) It is easy to see that
$$\align &n=p_3(x)+2p_4(y)+p_9(z)\ \t{for some}\ x,y,z\in\Z
\\\iff&56n+32=7(2x+1)^2+112y^2+(14z-5)^2\ \t{for some}\ x,y,z\in\Z
\\\iff&56n+32=7x^2+28y^2+z^2\ \ \t{for some}\ x,y,z\in\Z\ \t{with}\ 2\nmid z.
\endalign$$
(If $56n+32=7x^2+28y^2+z^2$ with $2\nmid z$, then  $28 y^2\eq32-7x^2-z^2\eq0\ (\mo\ 8)$ and hence $y$
is even, also $z^2\eq 32\eq2^2\ (\mo\ 7)$ and hence $z\eq \pm 5\ (\mo\ 14)$.)

Under Conjecture 1.15, one can check that all the numbers in $S(7)$ congruent to 32 mod 56 can be written in the form
$7x^2+28y^2+z^2$ with $z$ odd. If $56n+32\not\in S(7)$, then $56n+32=p+7x^2$ for some prime $p$ and $x\in\Z$.
Clearly $2\nmid x$ and $p\eq 32-7x^2\eq25\ (\mo\ 56)$. By (2.17) of [C, p.\,31],
there are $y,z\in\Z$ such that $p=7y^2+z^2$. Since $p\eq1\not\eq7\ (\mo\ 4)$, we have $2\mid y$.
So $56n+32=7x^2+28(y/2)^2+z^2$ with $2\nmid z$.

(ii) Observe that
$$\align &n=p_3(x)+p_4(y)+p_{13}(z)\ \t{for some}\ x,y,z\in\Z
\\\iff&88n+92=11(2x+1)^2+88y^2+(22z-9)^2\ \t{for some}\ x,y,z\in\Z
\\\iff&88(n+1)+4=11x^2+22y^2+z^2\ \ \t{for some}\ x,y,z\in\Z\ \t{with}\ 2\nmid z.
\endalign$$
(If $88(n+1)+4=11x^2+22y^2+z^2$ with $2\nmid z$, then  $22 y^2\eq4-11x^2-z^2\eq0\ (\mo\ 8)$ and hence $y$
is even, also $z^2\eq 4\ (\mo\ 11)$ and hence $z\eq \pm 9\ (\mo\ 22)$.)

Under Conjecture 1.15, one can check that all the numbers in $S(11)$ congruent to 4 mod 88 can be written in the form
$11x^2+22y^2+z^2$ with $z$ odd. If $88(n+1)+4\not\in S(11)$, then $88(n+1)+4=p+11x^2$
for some prime $p$ and $x\in\Z$.
Evidently, $2\nmid x$ and $p\eq 4-11x^2\eq 81\ (\mo\ 88)$. By (2.28) of [C, p.\,36],
there are $y,z\in\Z$ such that $p=22y^2+z^2$ and hence
$80(n+1)+4=11x^2+22y^2+z^2$ with $2\nmid z$.

In view of the above, we have completed the proof of Theorem 1.17.  \qed
\medskip

Recall the following results from [D27] and [D39, pp.\,112-113]:
$$\align E(5x^2+y^2+z^2)=&\{4^k(8l+3):\ k,l\in\N\},\tag5.1
\\E(24x^2+y^2+z^2)=&\bigcup_{k,l\in\N}\{4l+3,\,8l+6,\,9^k(9l+3)\},\tag5.2
\\E(3x^2+y^2+z^2)=&\{9^k(9l+6):\ k,l\in\N\},\tag5.3
\\E(16x^2+16y^2+z^2)=&\bigcup_{k,l\in\N}\{4l+2,4l+3,8l+5,16l+8,16l+12,4^k(8l+7)\}.\tag5.4
\endalign$$

\medskip
\noindent{\it Proof of Theorem 1.20}. By Theorem 1.9, none of the 6 sums in Theorem 1.20 is universal over $\N$.
Below we show that all the 6 sums are universal over $\Z$.
For convenience we fix a nonnegative integer $n$.

By (5.1), $40n+23=u^2+v^2+5w^2$ for some $u,v,w\in\Z$. The equality modulo $4$ yields that $2\nmid uvw$.
Note that $x^2\not\eq23\ (\mo\ 5)$ for any $x\in\Z$. So $u$ and $v$
are relatively prime to $10$. Since $u^2+v^2\eq23\eq-2\ (\mo\ 5)$, we
must have $u^2\eq v^2\eq-1\eq3^2\ (\mo\ 5)$. So $u$ or $-u$ has the
form $10y-3$ with $y\in\Z$, and $v$ or $-v$ has the form $10z-3$
with $z\in\Z$. Write $w=2x+1$ with $x\in\Z$. Then
$$40n+23=5(2x+1)^2+(10y-3)^2+(10z-3)^2$$
and hence $n=p_3(x)+p_7(y)+p_7(z).$

In light of (5.2), $24n+2=24x^2+u^2+v^2$ for some $u,v,x\in\Z$. Clearly $2\nmid uv$. Since $u^2+v^2\eq2\ (\mo\
3)$, we have $\gcd(uv,6)=1$. Thus, for some $y,z\in\Z$ we have
$$24n+2=24x^2+(6y-1)^2+(6z-1)^2$$
and hence $n=p_4(x)+p_5(y)+p_5(z).$

By (4.4), there are $w,x,y\in\Z$ such that $3n+1=w^2+3x^2+3y^2$. Note
that $w$ or $-w$ can be written as $3z-1$ with $z\in\Z$. So
$$3n+1=3x^2+3y^2+(3z-1)^2\ \t{and hence}\ n=p_4(x)+p_4(y)+p_8(z).$$

In view of (5.3), there are $u,v,x\in\Z$ such that $3n+2=3x^2+u^2+v^2$. As
$u^2+v^2\eq2\ (\mo\ 3)$, we have $\gcd(uv,3)=1$. So, for some $y,z\in\Z$ we have
$$3n+2=3x^2+(3y-1)^2+(3z-1)^2,\ \ \t{i.e.,}\ n=p_4(x)+p_8(y)+p_8(z).$$

By (5.4), $16n+9=w^2+16x^2+16y^2$ for some $w,x,y\in\Z$. As $w^2\eq3^2\pmod{16}$, $w$ or $-w$
has the form $8z-3$ with $z\in\Z$. Therefore
$$16n+9=16x^2+16y^2+(8z-3)^2\ \t{and hence}\ n=p_4(x)+p_4(y)+p_{10}(z).$$

By (4.10) there are $u,v,w\in\Z$ such that $48n+31=2u^2+2v^2+3w^2$.
The equality modulo $8$ yields that $2\nmid uvw$. Clearly $3\nmid uv$. Note that $3w^2\eq 31-2u^2-2v^2\eq 27\ (\mo\
16)$ and hence $w\eq\pm 3\ (\mo\ 8)$. So, for some $x,y,z\in\Z$ we have
$$48n+31=2(6x-1)^2+2(6y-1)^2+3(8z-3)^2$$
and hence $n=p_5(x)+p_5(y)+p_{10}(z).$

In view of the above, we have completed the proof of Theorem 1.20. \qed

\heading{6. Two auxiliary theorems}\endheading

\proclaim{Theorem 6.1} Let $a,b,c\in\Z^+$ with $a\ls b\ls c$,
and let $i,j,k\in\{3,4,\ldots\}$ with $\max\{i,j,k\}\gs5$.
Suppose that $(ap_i,bp_j,cp_k)$ is universal with $\min\{ai,bj,ck\}>5$. Then
$(ap_i,bp_j,cp_k)=(p_8,2p_3,3p_4)$.
\endproclaim
\Proof. We first claim that
$$(a,b,c)\in\{(1,1,3),(1,2,2),(1,2,3),(1,2,4)\}.$$
In fact, as $(ap_i,bp_j,cp_k)$ is universal and $ap_i(2),bp_j(2),cp_k(2)>5$, the set
$$S=\{0,a\}+\{0,b\}+\{0,c\}$$
contains $[0,5]=\{0,1,2,3,4,5\}$.
As $1,2\in S$, we have $a=1$ and $b\ls 2$. Since $b+2\in S$, $c$ cannot be greater than $b+2$.
Clearly $5\in S$ implies that $(a,b,c)\not=(1,1,1),(1,1,2)$. This proves the claim.
\medskip

{\it Case} 1. $(a,b,c)=(1,1,3)$.

Without loss of generality we assume that $i\ls j$. Note that $i=ai\gs6$.
Since $(ap_i,bp_j,cp_k)$ is universal and $6\not\in S$, we must have $i=\min\{ai,bj,ck\}=6$.
As $8\in N(p_6,p_j,3p_k)$, $p_j(3)=3(j-1)>8$ and $3p_k(2)=3k>8$, the set $\{0,1,6\}+\{0,1,j\}+\{0,3\}$
contains 8 and hence $j\in\{7,8\}$. It is easy to verify that
$12\not\in N(p_6,p_7,3p_3)$, $13\not\in N(p_6,p_8,3p_3)$,
$17\not\in N(p_6,p_7,3p_4)$ and $37\not\in N(p_6,p_8,3p_4)$.
As
$$12\not\in\{0,1,6\}+\{0,1,7\}+\{0,3\}\ \t{and}\ 13\not\in\{0,1,6\}+\{0,1,8\}+\{0,3\},$$
for $k\gs 5$ we have $12\not\in N(p_6,p_7,3p_k)$ and $13\not\in N(p_6,p_8,3p_k)$.
This contradicts the condition that
$(ap_i,bp_j,cp_k)$ is universal.
\medskip

{\it Case} 2. $(a,b,c)=(1,2,2)$.

Without loss of generality we assume that $j\ls k$.
Recall that $i,2j,2k\gs6$. Since $N(p_i,2p_j,2p_k)\supseteq[6,7]$, we have
$6,7\in\{0,1,i\}+\{0,2,2j\}+\{0,2,2k\}$.
As $\{0,1,i\}+\{0,2\}+\{0,2\}$ cannot contain both $6$ and $7$, we must have $2j=6$.
Observe that $p_i(3)=3(i-1)\gs15$ and $2p_k(3)\gs 2p_j(3)=6(j-1)=12$.
By $N(p_i,2p_j,2p_k)\supseteq[10,11]$, we get
$$10,11\in\{0,1,i\}+\{0,2,6\}+\{0,2,2k\}.$$
Since $\{0,1,i\}+\{0,2,6\}+\{0,2\}$ cannot contain both 10 and 11, we must have $2k<12$ and hence $k<6$.
Similarly, $k\not=3$ as $\{10,11\}\not\subseteq\{0,1,i\}+\{0,2,6\}+\{0,2,6\}$. So $k\in\{4,5\}$.

If $i\gs 17$, then $16\not\in N(p_i,2p_3,2p_4)$ since $\{0,1\}+\{0,2,6,12\}+\{0,2,8\}$ does not contain 16.
We can easily verify that $16\not\in N(p_i,2p_3,2p_4)$ for $i=7,9,11,13,15$ and that
$17\not\in N(p_i,2p_3,2p_4)$ for $i=8,10,12,14,16$. Also, $41\not\in N(p_6,2p_3,2p_4)$.
Similarly, $18\not\in N(p_i,2p_3,2p_5)$ for $i\gs 19$, since $18\not\in\{0,1\}+\{0,2,6,12\}+\{0,2,10\}$.
Also, $63\not\in N(p_6,2p_3,2p_5)$, $35\not\in N(p_7,2p_3,2p_5)$,
$19\not\in N(p_i,2p_3,2p_5)$ for
$i=8,10,12,14,16,18$, and $18\not\in N(p_i,2p_3,2p_5)$ for $i=9,11,13,15,17$.
Thus, when $k\in\{4,5\}$ we also have a contradiction.
\medskip

{\it Case} 3. $(a,b,c)=(1,2,3)$.

 Since $(p_i,2p_j,3p_k)$ is universal and $i,2j,3k\gs6$, we have
 $$7,8\in\{0,1,i\}+\{0,2,2j\}+\{0,3,3k\}.$$
 When $i=7$, we have $2j=8$ since $8\in\{0,1,7\}+\{0,2,2j\}+\{0,3\}$.
 Note that $13\not\in N(p_7,2p_4,3p_3)$ and $16\not\in N(p_7,2p_4,3p_4)$.
 If $k\gs 5$ then $13\not\in N(p_7,2p_4,3p_k)$ as $13\not\in\{0,1,7\}+\{0,2,8\}+\{0,3\}$.
 If $i\not=7$, then $7\not\in\{0,1,i\}+\{0,2\}+\{0,3\}$ and hence $2j=6$.
 As $8\not\in\{0,1\}+\{0,2,6\}+\{0,3\}$, we cannot have $i>8$.
 Thus $i\in\{6,8\}$ and $j=3$.

  Observe that $14\not\in N(p_6,2p_3,3p_3)$ and $22\not\in N(p_6,2p_3,3p_4)$.
  For $k\gs 5$ we have $14\not\in N(p_6,2p_3,3p_k)$ since $14\not\in\{0,1,6\}+\{0,2,6,12\}+\{0,3\}$.
 Note that
 $$35\not\in N(p_8,2p_3,3p_3),\ 19\not\in N(p_8,2p_3,3p_5),\ 22\not\in N(p_8,2p_3,3p_6).$$
 For $k\gs 7$ we have $18\not\in N(p_8,2p_3,3p_k)$ since $18\not\in \{0,1,8\}+\{0,2,6,12\}+\{0,3\}$.
  Thus $(i,j,k)=(8,3,4)$ and hence $(ap_i,bp_j,cp_k)=(p_8,2p_3,3p_4)$.
\medskip

{\it Case} 4. $(a,b,c)=(1,2,4)$.

 As $p_i(3)=3(i-1)\gs15$, $2p_j(3)=6(j-1)\gs 12$ and $4p_k(2)=4k\gs 12$, the set
 $$ T=\{0,1,i\}+\{0,2,2j\}+\{0,4\}$$
 must contain $8,9,10,11$. As $\{1,5\}+\{0,2,2j\}$ cannot contain both $9$ and 11,
 $i$ must be odd. Thus $8,10\in \{0,2,2j\}+\{0,4\}$, which is impossible.

 Combining the above we have completed the proof of Theorem 6.1. \qed

 \proclaim{Theorem 6.2} Let $b,c\in\Z^+$ with $b\ls c$, and let $j,k\in\{3,4,5,\ldots\}$ with $bj,ck\gs5$.
If $(p_5,bp_j,cp_k)$ is universal, then $(p_5,bp_j,cp_k)$ is on the following list:
$$(p_5,p_6,2p_4),\ (p_5,p_9,2p_3),\ (p_5,p_7,3p_3),\ (p_5,2p_3,3p_3),\ (p_5,2p_3,3p_4).$$
 \endproclaim
\Proof. As $(p_5,bp_j,cp_k)$ is universal, $bp_j(2)=bj\gs5$ and $cp_k(2)=ck\gs5$,  we have
$[0,4]\subseteq\{0,1\}+\{0,b,c,b+c\},$
hence $b\ls 2$ and $c\ls b+2$. Thus
$$(b,c)\in\{(1,2),(1,3),(2,2),(2,3),(2,4)\}.$$
Observe that $p_5(3)=3(5-1)=12$, $bp_j(3)=3b(j-1)\gs12$ and
$cp_k(3)=3c(k-1)\gs12$. Therefore
 $$S:=\{0,1,5\}+\{0,b,bj\}+\{0,c,ck\}\supseteq[6,11].$$
\medskip

{\it Case} 1. $(b,c)=(1,2)$.

 We first consider the case $2\mid j$.
 If $j=6$, then by $10\in S$ we get $2k=8$ and hence $(p_5,bp_j,cp_k)=(p_5,p_6,2p_4)$.
 For $k\gs9$ we have $16\not\in N(p_5,p_8,2p_k)$ since $16\not\in\{0,5,12\}+\{0,8\}+\{0,2\}$.
 For $k\in[3,8]$ it is easy to verify that $a_k\not\in N(p_5,p_8,2p_k)$, where
 $a_3=a_6=16$, $a_5=a_7=17$, $a_4=46$ and $a_8=19$.
 If $j=10$, then $9\in S$ implies that $k=4$. Note that $16\not\in N(p_5,p_{10},2p_4)$.
 For $j=12,14,16,\ldots$, the set $S$ cannot contain both $9$ and $11$.

 Now we consider the case $2\nmid j$. If $j=5$, then $S$ cannot contain both 9 and 11.
 If $j=7$, then $11\in S$ implies that $k\in\{3,5\}$. But $16\not\in N(p_5,p_7,2p_3)$ and $27\not\in N(p_5,p_7,2p_5)$.
 If $j=9$ and $k=3$ then $(p_5,bp_j,cp_k)=(p_5,p_9,2p_3)$. If $k\gs9$ then
 $17\not\in N(p_5,p_9,2p_k)$ since $17\not\in\{0,5,12\}+\{0,9\}+\{0,2\}$.
 For $k\in[4,8]$ it is easy to verify that $b_k\not\in N(p_5,p_9,2p_k)$, where
 $b_4=106$, $b_5=b_7=17$ and  $b_6=b_8=19$.
 For $j=11,13,15,\ldots$ we have $k=4$ by $9,10\in S$. Note that $17\not\in N(p_5,p_{11},2p_4)$
 and $11\not\in N(p_5,p_j,2p_4)$ for $j=13,15,17,\ldots$.
 \medskip

{\it Case} 2. $(b,c)=(1,3)$.

In this case, $ck\gs9$. By $7\in S$, we have $j\in\{6,7\}$.
It is easy to verify that
$$17\not\in N(p_5,p_6,3p_3),\ 65\not\in N(p_5,p_6,3p_4)\ \t{and}\ 24\not\in N(p_5,p_6,3p_5).$$
For $k\gs6$ we have $17\not\in N(p_5,p_6,3p_k)$ since
$$17\not\in\{0,1,5,12\}+\{0,1,6,15\}+\{0,3\}.$$
Note that $(p_5,p_7,3p_3)$ is on the list given in Theorem 6.2.
It is easy to check that $44\not\in N(p_5,p_7,3p_4)$. Also,
$14\not\in N(p_5,p_7,3p_k)$ for $k\gs5$, since $14\not\in\{0,1,5,12\}+\{0,1,7\}+\{0,3\}$.
\medskip

{\it Case} 3. $b=2$ and $c\in\{2,4\}$.

Since $b$ and $c$ are even, by $6,8,10\in S$ we get $3,4,5\in\{0,1,j\}+\{0,c/2,ck/2\}$.
This is impossible when $c=4$. Thus we let $c=2$ and assume $j\ls k$ without loss of generality.
By $3,4,5\in\{0,1,j\}+\{0,1,k\}$, we have $j=3$ and $k\in\{4,5\}$. One can verify that
$138\not\in N(p_5,2p_3,2p_4)$ and $60\not\in N(p_5,2p_3,2p_5)$.

\medskip

{\it Case} 4. $(b,c)=(2,3)$.

We first consider the case $k=3$. If $j=3$, then $(p_5,bp_j,cp_k)=(p_5,2p_3,3p_3)$.
Observe that
$$34\not\in N(p_5,2p_4,3p_3),\ 26\not\in N(p_5,2p_5,3p_3),\ 28\not\in N(p_5,2p_6,3p_3).$$
For $j\gs7$ we have $13\not\in N(p_5,2p_j,3p_3)$ since $13\not\in\{0,1,5,12\}+\{0,2\}+\{0,3,9\}$.

Now we consider the case $k\gs4$. Since $ck\gs12$, by $9\in S$ we get $2j\in\{6,8\}$.
Note that $(p_5,2p_3,3p_4)$ is on the list given in Theorem 6.2. Also, $19\not\in N(p_5,2p_3,3p_5)$
and $26\not\in N(p_5,2p_3,3p_6)$. For $k\gs7$ we have $19\not\in N(p_5,2p_3,3p_k)$ since
$$19\not\in\{0,1,5,12\}+\{0,2,6,12\}+\{0,3\}.$$
Note that $139\not\in N(p_5,2p_4,3p_4)$, $31\not\in N(p_5,2p_4,3p_k)$ for $k=5,7,9$,
and $28\not\in N(p_5,2p_4,3p_k)$ for $k=6,8$. For $k\gs10$ we have $28\not\in N(p_5,2p_4,3p_k)$ since
$$28\not\in\{0,1,5,12,22\}+\{0,2,8,18\}+\{0,3\}.$$

In view of the above we have finished the proof of Theorem 6.2. \qed

\heading{7. Proof of Theorem 1.9}\endheading

\proclaim{Lemma 7.1} Suppose that $(p_3,p_j,p_k)$ is universal with $3\ls j\ls k$ and $k\gs 5$.
Then $(j,k)$ is among the following ordered pairs:
$$\align &(3,k)\ (k=5,6,7,8,10,12,17),
\\&(4,k)\ (k=5,6,7,8,9,10,11,12,13,15,17,18,27),
\\&(5,k)\ (k=5,6,7,8,9,11,13),
\\&(7,8),\ (7,10).
\endalign$$
\endproclaim
\Proof. We distinguish four cases.
\medskip

{\it Case} 1. $j=3$.

 It is easy to verify that that $a_k\not\in N(p_3,p_3,p_k)$ for
 $$k=9,11,13,14,15,16,18,19,20,\ldots,33,$$ where
$$\align&a_9=a_{15}=a_{18}=a_{22}=a_{24}=a_{27}=a_{33}=41,\ a_{11}=a_{23}=63,
\\&a_{13}=a_{20}=a_{21}=a_{30}=53,\ a_{14}=a_{16}=a_{19}=a_{25}=a_{28}=33,
\\&a_{26}=129,\ a_{29}=125,\ a_{31}=54,\ a_{32}=86.
\endalign$$
For $k\gs 34$ we have $33\not\in N(p_3,p_3,p_k)$ since
$$33\not\in\{0,1,3,6,10,15,21,28\}+\{0,1,3,6,10,15,21,28\}+\{0,1\}.$$
\medskip

{\it Case} 2. $j=4$.

One can verify that that $b_k\not\in N(p_3,p_4,p_k)$ for
$$k=14,16,19,20,21,\ldots,26,28,29,\ldots,34,$$
where
$$\align&b_{14}=b_{16}=b_{21}=b_{26}=34,\ b_{19}=412,\ b_{20}=468,
\\&b_{22}=b_{32}=90,\ b_{23}=99,\ b_{24}=112,\ b_{25}=b_{28}=b_{30}=48,
\\&b_{29}=63,\ b_{31}=69,\ b_{33}=438,\ b_{34}=133.
\endalign$$
For $k\gs 35$ we have $34\not\in N(p_3,p_4,p_k)$ since
$$34\not\in\{0,1,3,6,10,15,21,28\}+\{0,1,4,9,16,25\}+\{0,1\}.$$
\medskip

{\it Case} 3. $j=5$.

 It is easy to verify that $c_k\not\in N(p_3,p_5,p_k)$ for
 $k=10,12,14,15,\ldots,31$, where
$$\align&c_{10}=c_{16}=c_{25}=c_{27}=c_{30}=69,\ c_{12}=c_{14}=c_{17}=c_{22}=31,
\\&c_{15}=c_{20}=131,\ c_{18}=c_{23}=c_{26}=c_{31}=65,\ c_{19}=168,
\\&c_{21}=135,\ c_{24}=218,\ c_{28}=75,\ c_{29}=82.
\endalign$$
For $k\gs 32$ we have $31\not\in N(p_3,p_5,p_k)$, since
$$31\not\in\{0,1,3,6,10,15,21,28\}+\{0,1,5,12,22\}+\{0,1\}.$$
\medskip

{\it Case} 4. $j\gs6$.

 Since $p_j(3)=3(j-1)\gs 15$ and $p_k(3)=3(k-1)\gs15$, we should have
 $$S=\{0,1,3,6,10\}+\{0,1,j\}+\{0,1,k\}\supseteq\{9,13,14\}.$$
 By $9\in S$ we get $j\ls 9$.

When $j=6$, by $14\in S$ we have $k\in\{7,8,10,11,12,13,14\}$.
Observe that $d_k\not\in N(p_3,p_6,p_k)$ for $k=7,8,10,11,12,13,14$, where
$$d_7=75,\ d_8=398,\ d_{10}=d_{11}=24,\ d_{12}=20,\ d_{13}=d_{14}=33.$$

Now we handle the case $j=7$. For $k=7,9,11,12,\ldots,26$ we have $e_k\not\in N(p_3,p_7,p_k)$,
where
$$\align&e_7=e_{13}=e_{15}=e_{18}=e_{22}=27,\ e_9=e_{24}=51,
\\&e_{11}=e_{16}=42,\ e_{12}=e_{14}=e_{17}=e_{21}=26,
\\&e_{19}=e_{26}=31,\ e_{20}=e_{23}=32,\ e_{25}=48.
\endalign$$
For $k\gs 27$ we have $26\not\in N(p_3,p_7,p_k)$ since
$$26\not\in\{0,1,3,6,10,15,21\}+\{0,1,7,18\}+\{0,1\}.$$

When $j=8$, by $13\in S$ we obtain $9\ls k\ls 13$. Observe that
$f_k\not\in N(p_3,p_8,p_k)$ for $k\in[9,13]$, where
$$f_9=413,\ f_{10}=104,\ f_{11}=84,\ f_{12}=59,\ f_{13}=26.$$

When $j=9$, by $14\in S$ we get $10\ls k\ls 14$. Note that
$g_k\not\in N(p_3,p_9,p_k)$ for $k\in[10,14]$, where
$g_{10}=g_{13}=18$, $g_{11}=43$ and $g_{12}=g_{14}=32.$

Combining the above we have proved the desired result. \qed

\proclaim{Lemma 7.2} Assume that $(p_4,p_j,p_k)$ is universal with $4\ls j\ls k$ and $k\gs 5$.
Then $j\in\{4,5\}$ and $k=j+1$.
\endproclaim
\Proof. As $p_j(3)=3(j-1)\gs 9$ and $p_k(3)=3(k-1)\gs9$, the set
$$S:=\{0,1,4\}+\{0,1,j\}+\{0,1,k\}$$
must contain $7$ and $8$. By $7\in S$ we see that $j\ls 7$.

When $j=4$, we have $k\in\{5,6,7\}$ by $7\in S$.
Note that $12\not\in N(p_4,p_4,p_6)$ and $77\not\in N(p_4,p_4,p_7)$.

When $j=5$, we have $k\in\{6,7,8\}$ by $8\in S$. Observe that
$63\not\in N(p_4,p_5,p_7)$ and $19\not\in N(p_4,p_5,p_8)$.

When $j=6$, one can verify that $s_k\not\in N(p_4,p_6,p_k)$ for $k\in[6,12]$, where
$$s_6=s_{11}=14,\ s_7=21,\ s_8=35,\ s_9=12,\ s_{10}=13,\ s_{12}=60.$$
For $k\gs 13$ we have $12\not\in N(p_4,p_6,p_k)$ since $12\not\in\{0,1,4,9\}+\{0,1,6\}+\{0,1\}$.

When $j=7$, it is easy to verify that $t_k\not\in N(p_4,p_7,p_k)$ for $k\in[7,13]$, where
$$t_7=t_{10}=13,\ t_8=t_{11}=14,\ t_9=t_{12}=15,\ t_{13}=42.$$
For $k\gs 14$ we have $13\not\in N(p_4,p_7,p_k)$ since $13\not\in\{0,1,4,9\}+\{0,1,7\}+\{0,1\}$.

In view of the above, we obtain that $(j,k)\in\{(4,5),(5,6)\}$.  \qed

\medskip

\noindent{\it Proof of Theorem 1.9}. By Theorems 6.1 and 6.2, we have $i\ls 5$ and $i\not=5$.
So $i\in\{3,4\}$. Thus the desired result follows from Lemmas 7.1 and 7.2. \qed

\heading{8. Proof of Theorem 1.11}\endheading

\proclaim{Lemma 8.1} Suppose that $(p_3,p_j,2p_k)$ is universal with $j,k\gs3$ and $\max\{j,k\}\gs 5$.
Then $(j,k)$ is among the following ordered pairs:
$$\align
&(j,3)\ (j=5,6,7,8,9,10,12,15,16,17,23),
\\&(5,4),\ (6,4),\ (8,4),\ (9,4),\ (17,4),
\\&(3,5),\ (4,5),\ (6,5),\ (7,5),\ (8,6),
\\&(4,7),\ (5,7),\  (7,7),\  (4,8),\ (4,9), \ (5,9).
\endalign$$
\endproclaim
\Proof. We distinguish six cases.
\medskip

{\it Case} 1. $k=3$.

In this case, $j\gs 5$. It is easy to see that
$a_j\not\in N(p_3,p_j,2p_3)$ for every $j=11,13,14,18,19,20,21,22,24,25$, where
$$\align&a_{11}=a_{14}=a_{21}=25,\ a_{13}=a_{18}=a_{25}=50,
\\&a_{19}=258,\ a_{20}=89,\ a_{22}=54,\ a_{24}=175.
\endalign$$
For $j\gs 26$ we have $25\not\in N(p_3,p_j,2p_3)$ since
$$25\not\in\{0,1,3,6,10,15,21\}+\{0,1\}+\{0,2,6,12,20\}.$$
\medskip

{\it Case} 2. $k=4$.

 One can verify that $b_j\not\in N(p_3,p_j,2p_4)$ for
 $$j=7,\,10,\,11,\ldots,16,\,18,\,19,\ldots,26,$$ where
$$\align&b_7=59,\ b_{10}=b_{13}=b_{19}=b_{22}=26,\ b_{11}=b_{14}=b_{20}=b_{23}=27,
\\&b_{12}=b_{18}=b_{24}=49,\ b_{15}=b_{16}=b_{21}=b_{25}=41,\ b_{26}=115.
\endalign$$
For $j\gs 27$ we have $26\not\in N(p_3,p_j,2p_4)$ since
$$26\not\in\{0,1,3,6,10,15,21\}+\{0,1\}+\{0,2,8,18\}.$$
\medskip

{\it Case} 3. $k=5$.

Observe that $c_j\not\in N(p_3,p_j,2p_5)$ for $j=5,8,9,\ldots,19$, where
$$\align&c_5=c_{10}=c_{12}=c_{15}=19,\ c_8=83,\ c_9=c_{16}=42,
\\&c_{11}=c_{13}=62,\ c_{14}=c_{19}=33,\ c_{17}=43,\ c_{18}=114.
\endalign$$
For $j\gs20$ we have $19\not\in N(p_3,p_j,2p_5)$ since
$$19\not\in\{0,1,3,6,10,15\}+\{0,1\}+\{0,2,10\}.$$
\medskip

{\it Case} 4. $k=6$.

It is easy to see that $d_j\not\in N(p_3,p_j,2p_6)$ for $j=3,4,5,6,7,9,
10,\ldots,20$, where
$$\align&d_3=35,\ d_4=d_9=d_{11}=d_{13}=d_{16}=20,
\\&d_{5}=124,\ d_6=d_{10}=d_{12}=d_{15}=d_{17}=d_{19}=26,
\\&d_7=50,\ d_{14}=d_{18}=25,\ d_{20}=44.
\endalign$$
For $j\gs21$ we have $20\not\in N(p_3,p_j,2p_6)$ since
$$20\not\in\{0,1,3,6,10,15\}+\{0,1\}+\{0,2,12\}.$$
\medskip

{\it Case} 5. $k=7$.

Observe that $e_j\not\in N(p_3,p_j,2p_7)$ for $j=3,6,8,9,\ldots,19$, where
$$\align&e_3=e_6=e_8=e_{10}=e_{12}=e_{15}=19,\ e_9=86,
\\&e_{11}=e_{14}=e_{16}=e_{18}=27,
\ e_{13}=e_{17}=e_{19}=26.
\endalign$$
For $j\gs20$ we have $19\not\in N(p_3,p_j,2p_7)$ since
$$19\not\in\{0,1,3,6,10,15\}+\{0,1\}+\{0,2,14\}.$$
\medskip

{\it Case} 6. $k\gs8$.

As $14\in N(p_3,p_j,2p_k)$, $p_j(5)=10j-15\gs15$ and $2k>14$, we have
$$14\in\{0,1,3,6,10\}+\{0,1,j,3j-3,6j-8\}+\{0,2\}.$$
It follows that $j\in\{3,4,5,6,8,9,11,12,13,14\}$.

For $j=3,5,6,8,12$ it is easy to see that $19\not\in N(p_3,p_j,2p_k)$ for $k\gs 10$.
Note that $35\not\in N(p_3,p_j,2p_k)$ for $j\in\{3,8\}$ and $k\in\{8,9\}$.
Also, $26\not\in N(p_3,p_j,2p_9)$ for $j=6,12$, and
$$124\not\in N(p_3,p_5,2p_8),\ 35\not\in N(p_3,p_6,2p_8),\ 25\not\in N(p_3,p_{12},2p_8).$$

For $j=4,9,11,13$ it is easy to see that $20\not\in N(p_3,p_j,2p_k)$ for $k\gs 11$.
Observe that $43\not\in N(p_3,p_4,2p_{10})$. Also,
$$\align&33\not\in N(p_3,p_9,2p_8),\ 35\not\in N(p_3,p_9,2p_9),\ 33\not\in N(p_3,p_9,2p_{10}),
\\&25\not\in N(p_3,p_{11},2p_8),\ 27\not\in N(p_3,p_{11},2p_9),\ 25\not\in N(p_3,p_{11},2p_{10}).
\\&33\not\in N(p_3,p_{13},2p_8),\ 26\not\in N(p_3,p_{13},2p_9),\ 32\not\in N(p_3,p_{13},2p_{10}).
\endalign$$

For the case $j=14$, we have
$$27\not\in N(p_3,p_{14},2p_9),\ 27\not\in N(p_3,p_{14},2p_{11}),\ 32\not\in N(p_3,p_{14},2p_{12}).$$
Also, $25\not\in N(p_3,p_{14},2p_k)$ for $k=8,10,13,14,\ldots$.

\medskip

Combining the above we have completed the proof. \qed

\proclaim{Lemma 8.2} Let $c$ be a positive integer greater than $2$.
Suppose that $(p_3,p_j,cp_k)$ is universal with $j,k\gs3$ and $\max\{j,k\}\gs 5$.
Then $(p_j,cp_k)$ is among the following $10$ ordered pairs:
$$\align
&(p_3,4p_5)\ (p_5,3p_3),\ (p_5,4p_3),\ (p_5,6p_3),\ (p_5,9p_3),
\\&(p_5,3p_4),\ (p_5,4p_4),\ (p_5,4p_6),\ (p_5,4p_7),\ (p_8,3p_4).
\endalign$$
\endproclaim
\Proof. Since $(p_3,p_j,cp_k)$ is universal, $cp_k(3)=3c(k-1)\gs 9(k-1)\gs18$ and $p_j(6)=15j-24\gs21$,
the set
$$T=\{0,1,3,6,10,15\}+\{0,1,j,3j-3,6j-8,10j-15\}+\{0,c,ck\}$$
must contain $\{5,8,9,12,13,14,15,16,17\}$.

\medskip

{\it Case} 1. $c=3$.

By $8\in T$ we get $j\in\{4,5,7,8\}$.
Observe that $41\not\in N(p_3,p_4,3p_k)$ for $k=6,7$, and $23\not\in N(p_3,p_4,3p_k)$ for $k=5,8,9,\ldots$.
For $j=5$ we have $k\in\{3,4,5\}$ by $17\in T$. Note that $34\not\in N(p_3,p_5,3p_5)$.
For $j=7$, by $12\in T$ we obtain $k\in\{3,4\}$. It is easy to see that
$23\not\in N(p_3,p_7,3p_3)$ and $26\not\in N(p_3,p_7,3p_4)$.
Also, $35\not\in N(p_3,p_8,3p_k)$ for $k=3,6$, and $20\not\in N(p_3,p_8,3p_k)$ for $k=5,7,8,\ldots$.

\medskip

{\it Case} 2. $c=4$.

By $9\in T$ we have $j\ls 9$ and $j\not=7$. For $j\in\{6,8,9\}$, we have $k=4$ by $17\in T$.
Note that
$24\not\in N(p_3,p_6,4p_4)$, $98\not\in N(p_3,p_8,4p_4)$ and $84\not\in N(p_3,p_9,4p_4)$.
For $j=3$, we have $23\not\in N(p_3,p_3,4p_k)$ for $k\gs6$.
For $j=4$, we have $38\not\in N(p_3,p_4,4p_5)$, $47\not\in N(p_3,p_4,4p_6)$, and $27\not\in N(p_3,p_4,4p_k)$ for $k\gs7$.
For $j=5$, we have $143\not\in N(p_3,p_5,4p_8)$,
and $34\not\in N(p_3,p_5,4p_k)$ for $k=5,9,10,\ldots$.

\medskip

{\it Case} 3. $c=5$.

Since $5k\gs15$, by $13,14\in T$ we get $j\in\{3,8,13\}$.
If $j\in\{8,13\}$, then by $17\in T$ we obtain $k=3$.
Note that $35\not\in N(p_3,p_j,5p_3)$ for $j=8,13$.
Also, $19\not\in N(p_3,p_3,5p_k)$ for $k\gs5$.

\medskip

{\it Case} 4. $c>5$.

By $5\in T$ we have $j\in\{4,5\}$. If $j=4$ then $c\ls 8$ by $8\in T$.
If $j=5$, then $c\in\{6,7,9\}$ by $9,17\in T$. Observe that
$68\not\in N(p_3,p_4,6p_5)$, and $33\not\in N(p_3,p_4,6p_k)$ for $k\gs 6$.
Also, $68\not\in N(p_3,p_5,6p_4)$, $114\not\in N(p_3,p_5,6p_5)$,
and $30\not\in N(p_3,p_5,6p_k)$ for $k\gs 6$. Note that
$20\not\in N(p_3,p_4,7p_k)$ for $k\gs5$. Also, $89\not\in N(p_3,p_5,7p_3)$, and
$24\not\in N(p_3,p_5,7p_k)$ for $k\gs 4$. It is easy to verify that
$273\not\in N(p_3,p_4,8p_5)$ and $41\not\in N(p_3,p_4,8p_k)$ for $k\gs6$.
Also, $53\not\in N(p_3,p_5,9p_4)$ and $39\not\in N(p_3,p_5,9p_k)$ for $k\gs 5$.

\medskip

In view of the above, we obtain the desired result. \qed

\proclaim{Lemma 8.3} Let $b$ and $c$ be integers with $2\ls b\ls c$.
Let $j,k\in\{3,4,\ldots\}$, $\max\{j,k\}\gs 5$, and $j\ls k$ if $b=c$.
Suppose that $(p_3,bp_j,cp_k)$ is universal.
Then $(bp_j,cp_k)$ is among the following $10$ ordered pairs:
$$\align
&(2p_3,2p_6),\ (2p_3,2p_7),\ (2p_3,2p_8),\ (2p_3,2p_9),\ (2p_3,2p_{12}),
\\&(2p_3,4p_5),\ (2p_4,2p_5),\ (2p_4,4p_5), (2p_5,4p_3),\ (2p_5,4p_4).
\endalign$$
\endproclaim
\Proof. Since $(p_3,bp_j,cp_k)$ is universal and $bj,ck\gs2\times3=6$, we have
$2,4\in\{0,1,3\}+\{0,b\}+\{0,c\}$, which implies that $b=2$ and $c\in\{2,3,4\}$.

Suppose that $ck<12$. Then, either  $c=2$ and $j\ls k=5$, or $c=k=3$. Note that $139\not\in N(p_3,2p_3,2p_5)$,
$9\not\in N(p_3,2p_5,2p_5)$, and $7\not\in N(p_3,2p_j,3p_3)$ for $j\gs 5$.

Below we assume that $ck\gs 12$. As $(p_3,2p_j,cp_k)$ is universal and $2p_j(3)=6(j-1)\gs12$,  the set
$$R=\{0,1,3,6,10\}+\{0,2,2j\}+\{0,c\}$$
must contain $[0,11]$.

When $c=2$, we have $k\gs 6$ by $ck\gs12$, and $j\in\{3,4\}$ by $9\in R$.
Note that $r_k\not\in N(p_3,2p_3,2p_k)$ for $k=10,11,13,14,\ldots$, where
$$\align&r_{10}=r_{14}=r_{20}=r_{21}=\cdots=39,\ r_{11}=r_{16}=46,
\\&r_{13}=r_{15}=76,\ r_{17}=83,\ r_{18}=151,\ r_{19}=207.
\endalign$$
Also, $43\not\in N(p_3,2p_4,2p_6)$, $27\not\in N(p_3,2p_4,2p_k)$ for $k\in\{7,10\}$,
$64\not\in N(p_3,2p_4,2p_8)$, $826\not\in N(p_3,2p_4,2p_{11})$,
and $22\not\in N(p_3,2p_4,2p_k)$ for $k=9,12,13,\ldots$.

In the case $c=3$, by $7\in R$ we get $j=3$. Note that $14\not\in N(p_3,2p_3,3p_k)$ for $k\gs 5$.

When $c=4$, we have $j\in\{3,4,5\}$ by $11\in R$. Observe that $53\not\in N(p_3,2p_3,4p_k)$ for $k=6,7$,
and $29\not\in N(p_3,2p_3,4p_k)$ for $k\gs 8$. Also, $20\not\in N(p_3,2p_4,4p_k)$ for $k\gs6$,
and $18\not\in N(p_3,2p_5,4p_k)$ for $k\gs 5$.

By the above, we have completed the proof. \qed

\proclaim{Lemma 8.4} Let $b$ and $c$ be positive integers with $b\ls c$ and $c>1$.
Let $j,k\in\{3,4,5,\ldots\}$, $\max\{j,k\}\gs 5$, and $j\ls k$ if $b=c$.
Suppose that $(p_4,bp_j,cp_k)$ is universal with $bj\gs4$.
Then $(bp_j,cp_k)$ is among the following $11$ ordered pairs:
$$(p_j,2p_3)\,(j=5,6,7,8,10,12,17),\ (p_5,2p_4),\ (p_5,3p_3),\ (2p_3,2p_5),\ (2p_3,4p_5).
$$
\endproclaim
\Proof. Since $(p_4,bp_j,cp_k)$ is universal and $c$ is greater than one, we must have $2,3\in\{0,1\}+\{0,b\}+\{0,c\}$.
Thus, $b=1 $ and $c\in\{2,3\}$, or $b=2\ls c$.
\medskip

{\it Case} 1. $b=1$ and $c=2$.

 In view of (1.5),
$$\align&(p_4,p_j,2p_3)\ \t{is universal}
\\\iff&(p_3,p_3,p_j)\ \t{is universal}
\\\Longrightarrow\ &j\in\{5,6,7,8,10,12,17\}\ (\t{by Lemma 7.1}).
\endalign$$

Now let $k\gs 4$. Note that $j=bj\gs4$.

In the case $j=4$, we have
$21\not\in N(p_4,p_4,2p_k)$ for $k\in\{5,7\}$, $23\not\in N(p_4,p_4,2p_6)$,
and $14\not\in N(p_4,p_4,2p_k)$ for $k\gs 8$.

When $j=5$, we have $42\not\in N(p_4,p_5,2p_5)$, $29\not\in N(p_4,p_5,2p_k)$ for $k\in\{7,9\}$,
$34\not\in N(p_4,p_5,2p_8)$, $111\not\in N(p_4,p_5,2p_{10})$, and
$20\not\in N(p_4,p_5,2p_k)$ for $k=6,11,12,\ldots$.

For the case $j=6$, it is easy to verify that $80\not\in N(p_4,p_6,2p_4)$,
$20\not\in N(p_4,p_6,2p_6)$, and $13\not\in N(p_4,p_6,2p_k)$ for $k=5,7,8,\ldots$.

When $j=7$, we have $30\not\in N(p_4,p_7,2p_5)$, $15\not\in N(p_4,p_7,2p_6)$,
$42\not\in N(p_4,p_7,2p_7)$, and $14\not\in N(p_4,p_7,2p_k)$ for $k=4,8,9,\ldots$.

In the case $j=8$, we have $15\not\in N(p_4,p_8,2p_k)$ for $k=4,6$, and
$13\not\in N(p_4,p_8,2p_k)$ for $k=5,7,8,\ldots$.

When $j>8$, we have $k=4$ by $8\in N(p_4,p_j,2p_k)$. Note that
$n_j\not\in N(p_4,p_j,2p_4)$ for $j=9,10,\ldots$, where
$$n_9=n_{15}=n_{16}=\cdots=14,\ n_{10}=15,\ n_{11}=n_{14}=21,\ n_{12}=40,\ n_{13}=91.$$

\medskip

{\it Case} 2. $b=1$ and $c=3$.

Note that $p_3(j)=3(j-1)\gs9$ since $j=bj\gs4$. By $6\in N(p_4,p_j,3p_k)$, we have $6\in\{0,1,4\}+\{0,1,j\}+\{0,3\}$ and hence $j\in\{5,6\}$.
It is easy to verify that $11\not\in N(p_4,p_j,3p_k)$ for $j\in\{5,6\}$ and $k\gs 4$. Note also that
$21\not\in N(p_4,p_6,3p_3)$.

\medskip

{\it Case} 3. $b=2\ls c$.

Observe that
$$\align&(p_4,2p_3,cp_k)\ \t{is universal}
\\\iff&(p_3,p_3,cp_k)\ \t{is universal}\ \ (\t{by (1.5)})
\\\Longrightarrow\ &k=5\ \t{and}\ c\in\{2,4\}\ \ (\t{by Lemmas 8.1 and 8.2}).
\endalign$$

Now let $j\gs 4$. Clearly $2p_j(3)=6(j-1)\gs18$ and
$$cp_k(3)=3c(k-1)\gs\min\{6(j-1),9(k-1)\}\gs18.$$
Thus, by $[0,15]\subseteq N(p_4,2p_j,cp_k)$, the set
$$S=\{0,1,4,9\}+\{0,2,2j\}+\{0,c,ck\}$$
contains $[0,15]$. Note that $c\ls 5$ by $5\in S$.
If $c=2$, then $k\gs j\gs4$ and hence $7\not\in S$.
When $c=3$, we have $j=4$ by $8\in S$, hence $10\not\in S$ since $ck\gs15$.
If $c=4$ and $k\gs4$, then $\{12,14\}\not\subseteq S$.
Also, $19\not\in N(p_4,2p_j,4p_3)$ for $j=6,8$, and $17\not\in N(p_4,2p_j,4p_3)$ for $j=5,7,9,10,\ldots$.
When $c=5$, we have $\{10,12\}\not\subseteq S$.

In view of the above, we have proved Lemma 8.4. \qed

\medskip
\noindent{\it Proof of Theorem 1.11}.
 Combining Theorems 6.1-6.2 and Lemmas 8.1-8.4 we immediately obtain the desired result. \qed
\medskip

\Ack. The work was supported by the National Natural Science Foundation (Grant No.
11171140) of China and the PAPD of Jiangsu Higher Education Institutions.
The author thanks the two referees for their helpful comments.
The initial version of this paper was posted to arXiv as a preprint in May 2009 with the ID {\tt arXiv:0905.0635}.


\widestnumber\key{OnoS}

\Refs

\ref\key B\by B. C. Berndt\book  Number Theory in the Spirit of Ramanujan
\publ Amer. Math. Soc., Providence, R.I., 2006\endref

\ref\key CH\by W. K. Chan and A. Haensch\paper Almost universal ternary sums of squares and triangular numbers\jour in: Quadratic and Higher Degree Forms, in: Dev. Math., vol. 31, Springer, New York, 2013, pp. 51--62\endref

\ref\key CO\by W. K. Chan and B.-K. Oh\paper Almost universal ternary sums of triangular numbers\jour Proc. Amer. Math. Soc.\vol 137\yr 2009\pages 3553--3562\endref

\ref\key C\by D. A. Cox\book Primes of the Form $x^2+ny^2$\publ John Wiley \& Sons, New York, 1989\endref

\ref\key D27\by L. E. Dickson\paper Integers represented by positive ternary quadratic forms
\jour Bull. Amer. Math. Soc.\vol 33\yr 1927\pages 63--77\endref

\ref\key D39\by L. E. Dickson\book
Modern Elementary Theory of Numbers
\publ University of Chicago Press, Chicago, 1939\endref

\ref\key D99a\by L. E. Dickson\book
History of the Theory of Numbers, {\rm Vol. I}
\publ AMS Chelsea Publ., 1999\endref

\ref\key D99b\by L. E. Dickson\book
History of the Theory of Numbers, {\rm Vol. II}
\publ AMS Chelsea Publ., 1999\endref

\ref\key G\by E. Grosswald\book Representation of Integers as Sums of Squares
\publ Springer, New York, 1985\endref

\ref\key GPS\by S. Guo, H. Pan and Z.-W. Sun\paper
Mixed sums of squares and triangular numbers (II)
\jour Integers\vol 7\yr 2007\pages \#A56, 5pp (electronic)\endref

\ref\key Gu\by R. K. Guy\jour {\it Every number is expressible as the sum of how many polygonal numbers?}
Amer. Math. Monthly\vol 101\yr 1994\pages 169--172\endref

\ref\key JKS\by W. C. Jagy, I. Kaplansky and A. Schiemann\paper There are 913 regular ternary forms
\jour Mathematika\yr 1997\vol 44\pages 332--341\endref

\ref\key JP\by B. W. Jones and G. Pall
\paper Regular and semi-regular positive ternary quadratic forms
\jour Acta Math.\vol 70\yr 1939\pages 165--191\endref

\ref\key K09\by B. Kane\paper On two conjectures about mixed sums of squares and triangular numbers
\jour J. Comb. Number Theory\vol 1\yr 2009\pages 77--90\endref

\ref\key KS\by B. Kane and Z.-W. Sun\paper On almost universal mixed sums of squares and triangular numbers
\jour Trans. Amer. Math. Soc.\vol 362\yr 2010\pages 6425--6455\endref

\ref\key MW\by C. J. Moreno and S. S. Wagstaff\book Sums of Squares of Integers
\publ Chapman \& Hall/CRC, New York, 2005\endref

\ref\key N87\by M. B. Nathanson\paper A short proof of Cauchy's polygonal theorem
\jour Proc. Amer. Math. Soc.\vol 99\yr 1987\pages 22--24\endref

\ref\key N96\by M. B. Nathanson\paper Additive Number Theory: The
Classical Bases \publ Grad. Texts in Math., vol. 164, Springer,
New York, 1996\endref

\ref\key OS\by B.-K. Oh and Z.-W. Sun\paper Mixed sums of squares and triangular numbers (III)
\jour J. Number Theory\vol 129\yr 2009\pages 964-969\endref

\ref\key{OnoS}\by K. Ono and K. Soundararajan\paper
Ramanujan's ternary quadratic form\jour Invent. Math. \vol 130\yr 1997\pages 415--454.\endref

\ref\key S07\by Z.-W. Sun\paper Mixed sums of
squares and triangular numbers \jour Acta Arith. \vol 127\yr 2007\pages 103--113\endref

\ref\key S09\by Z.-W. Sun\paper On sums of primes and triangular
numbers \jour J. Comb. Number Theory \vol 1\yr
2009\pages 65--76\endref

\ref\key WP\by X. Wang and D. Pei\paper Eisenstein series of $3/2$ weight and one conjecture of Kaplansky
\jour Sci. China Ser. A \vol 44\yr 2001\pages 1278--1283\endref


\endRefs
\enddocument